\documentclass[letterpaper,11pt]{article}

\usepackage{graphicx}

\usepackage{amsmath,amssymb,relsize}
\usepackage{enumitem}
\usepackage{algorithm,algorithmic}
\usepackage{hyperref}
\usepackage{caption,subcaption}
\usepackage[titletoc,title]{appendix}
\usepackage{geometry} 
\usepackage{comment}

\AtBeginDocument{
  \newgeometry{
    textheight=9in,
    textwidth=6.5in,
    top=1in,
    headheight=14pt,
    headsep=25pt,
    footskip=30pt
  }
}

\def\grad{\nabla}

\def\bb{\mathbf{b}}

\def\br{\mathbf{r}}
\def\bs{\mathbf{s}}

\def\bu{\mathbf{u}}
\def\bv{\mathbf{v}}
\def\bw{\mathbf{w}}
\def\bx{\mathbf{x}}  
\def\by{\mathbf{y}}
\def\bz{\mathbf{z}}

\def\bo{\mathbf{0}}

\def\cB{\mathcal{B}}
\def\cC{\mathcal{C}}

\def\cF{\mathcal{F}}

\def\cN{\mathcal{N}}
\def\cO{\mathcal{O}}

\def\cS{\mathcal{S}}

\def\cX{\mathcal{X}}
\def\cY{\mathcal{Y}}

\def\mE{\mathbb{E}}

\def\mN{\mathbb{N}}

\def\mR{\mathbb{R}}

\def\smskip{\smallskip}

\def\texitem#1{\par\smskip\noindent\hangindent 25pt
               \hbox to 25pt {\hss #1 ~}\ignorespaces}

\def\abs#1{\left|#1\right|}
\def\norm#1{\left\|#1\right\|}

\newcommand{\BEAS}{\begin{eqnarray*}}
\newcommand{\EEAS}{\end{eqnarray*}}
\newcommand{\BEA}{\begin{eqnarray}}
\newcommand{\EEA}{\end{eqnarray}}
\newcommand{\BEQ}{\begin{eqnarray}}
\newcommand{\EEQ}{\end{eqnarray}}
\newcommand{\BIT}{\begin{itemize}}
\newcommand{\EIT}{\end{itemize}}
\newcommand{\BNUM}{\begin{enumerate}}
\newcommand{\ENUM}{\end{enumerate}}

\newcommand{\BA}{\begin{array}}
\newcommand{\EA}{\end{array}}

\def\fprod#1{\left\langle#1\right\rangle}

\newcommand\numberthis{\addtocounter{equation}{1}\tag{\theequation}}

\newcommand{\dist}{\mathop{\bf dist{}}}

\newcommand{\argmin}{\mathop{\rm argmin}}

\newcommand{\dom}{\mathop{\bf dom}}

\newif\ifpagenumbering
\pagenumberingtrue

\pagenumberingfalse

\newsavebox{\theorembox}
\newsavebox{\lemmabox}
\newsavebox{\defnbox}
\newsavebox{\assbox}
\savebox{\theorembox}{\noindent\bf Theorem}
\savebox{\lemmabox}{\noindent\bf Lemma}
\savebox{\defnbox}{\noindent\bf Definition}
\savebox{\assbox}{\noindent\bf Assumption}
\newtheorem{assumption}{\usebox{\assbox}}

\newcommand*{\affaddr}[1]{#1} 

\newcommand*{\email}[1]{\texttt{#1}}
\newtheorem{example}{Example}
\newtheorem{theorem}{Theorem}[section]
\newtheorem{corollary}{Corollary}[section]
\newtheorem{lemma}{Lemma}[section]
\newtheorem{remark}{Remark}[section]
\newtheorem{proposition}{Proposition}[section]
\newtheorem{definition}{Definition}

\newcommand{\qed}{\hfill \ensuremath{\Box}}

\usepackage{subfiles}

\begin{document}

\title{Stochastic Composition Optimization of Functions without Lipschitz Continuous Gradient}
\author{%
	Yin Liu and Sam Davanloo Tajbakhsh \\
	\email{\{liu.6630,davanloo.1\}@osu.edu}\\
	\affaddr{The Ohio State University}%
}

\maketitle

\begin{abstract}
	In this paper, we study stochastic optimization of two-level composition of functions without Lipschitz continuous gradient. The smoothness property is generalized by the notion of \textit{relative smoothness} which provokes the Bregman gradient method. We propose three Stochastic Composition Bregman Gradient algorithms for the three possible  {relatively smooth} compositional scenarios and provide their sample complexities to achieve an $\epsilon$-approximate stationary point. For the smooth of relatively smooth composition, the first algorithm requires $\cO(\epsilon^{-2})$ calls to the stochastic oracles of the inner function value and gradient as well as the outer function gradient.  When both functions are relatively smooth, the second algorithm requires $\cO(\epsilon^{-3})$ calls to the inner function value stochastic oracle and $\cO(\epsilon^{-2})$ calls to the inner and outer functions gradients stochastic oracles.  We further improve the second algorithm by variance reduction for the setting where just the inner function is smooth. The resulting algorithm requires $\cO(\epsilon^{-5/2})$ calls to the inner function value stochastic oracle, $\cO(\epsilon^{-3/2})$ calls to the inner function gradient and $\cO(\epsilon^{-2})$ calls to the outer function gradient stochastic oracles. Finally, we numerically evaluate the performance of these  {three} algorithms over two  {different} examples.
\end{abstract}

\vspace{0.2cm}

\textbf{keywords:} composition optimization, stochastic optimization algorithms, Bregman subproblem



\section{Introduction}
 This paper considers the two-level stochastic composition problem
\begin{align} \label{eq:2-level composition}
	\min_{\bx\in \cX} F(\bx) & \triangleq  f(g(\bx))
\quad \text{with } f(\bu)\triangleq \mathbb{E}[f_\varphi(\bu)], \ 
 g(\bx)\triangleq \mathbb{E}[g_\xi(\bx)],
\end{align}
where $\mathbb{E}[f_\varphi(\bu)]\triangleq\int_{\Omega_f}f_{\varphi(\omega_f)}dP_f(\omega_f)$, $\mathbb{E}[g_\xi(\bx)] \triangleq\int_{\Omega_g}g_{\xi(\omega_g)}dP_g(\omega_g)$, and $\cX$ is a closed convex set. Here $P_f$ (\emph{similarly} $P_g$) is a probability distribution on the sample space $\Omega_f$ (\emph{similarly} $\Omega_g$), the random vector $\varphi$ (\emph{similarly} $\xi$) is a mapping from $\Omega_f$ (\emph{similarly} $\Omega_g$) to a measurable space $\mathbb{W}_f$ (\emph{similarly} $\mathbb{W}_g$), and $f_\varphi\colon\mathbb{R}^d\rightarrow\mathbb{R}$ is a smooth function (\emph{similarly} $g_\xi\colon\mathbb{R}^n\rightarrow\mathbb{R}^d$ is a smooth map).

Furthermore, we assume that the functions $f$, $g$, or both do \emph{not} have Lipschitz continuous gradient. A continuously differentiable function $h$ has Lipschitz continuous gradient if for some constant $L$, we have 
\begin{align}\label{eq:lip_grad_def}
\norm{\grad h(\bx) - \grad h(\bar{\bx})} \leq L\norm{\bx - \bar{\bx}}, \quad \forall \bx,\ \bar{\bx} \in \dom h.
\end{align}

In the absence of Lipschitz continuity of the inner, outer, or both functions' gradients, problem \eqref{eq:2-level composition} appears in many applications such as policy evaluation for Markov decision processes~\cite{dann2014policy}, risk-averse optimization~\cite{ruszczynski2009risk}, low-rank non-negative matrix factorization~\cite{berry2007algorithms,ge2015escaping}, and model-agnostic meta-learning (MAML)~\cite{finn2017model}. Below, we discuss some of these applications in more detail.

\begin{example}[Policy Evaluation  for Markov decision processes (MDP)]\label{ex:MDP}
	Consider a Markov chain with states $\{Y_0,Y_1,\dots\}\subset \cY$, an unknown transition operator $P$, a reward function $r: \cY \rightarrow \mathbb{R}$ and a discount factor $\gamma \in (0,1)$. The goal is to estimate the value function $V: \cY \rightarrow\mathbb{R}$, defined as $V(y) \triangleq \mathbb{E} \left[ \sum_{t=0}^\infty \gamma^t r(Y_t) |Y_0 = y  \right]$, for a fixed control policy $\pi$. For a finite state space $\cY$, the value functions of all possible initial states can be represented as a vector $\bv\in \mathbb{R}^{|\cY|}$ that satisfies the Bellman equation $ \bv = \br+\gamma P\bv$,
	where $\br\in\mR^{|\cY|}$ such that $r_i=r(Y_i)$.
	When $|\cY|$ is large, one can approximate $\bv$ by $\bv \approx\Phi \bx$, where $\Phi$ is some matrix of basis functions and $\bx\in\mR^n$ contains the coefficients with $n\ll |\cY|$.
	Since one only has access to random samples of $P$ and $r$, i.e., $\hat{P}$ and $\hat{r}$, the residual minimization problem  to approximate the solution of the above system is
	\begin{align}\label{eq:ex_policy_eval}
	\min_{\bx\in\cX\subseteq\mR^n} \dist (\mathbb{E}[\hat{\br}],(I-\gamma \mathbb{E}[\hat{P}])\Phi\bx),
	\end{align}
	where $\dist:\mR^{|\cY|}\times\mR^{|\cY|}\rightarrow\mR_+$ is some distance function~\cite{sutton2018reinforcement}. Let $A_\xi \triangleq (I-\gamma\hat{P})\Phi$ and $\br_\varphi \triangleq \hat{\br} $, the problem can then be written as $\min_{\bx\in\cX} \dist (\mathbb{E}[\br_\varphi],\mathbb{E}[A_\xi]\bx)$ which is a special case of \eqref{eq:2-level composition}. Note that if {the} reward is in form of count data (e.g., number of news clicks suggested by a recommender system -- see, e.g., \cite{luo2022mindsim} and references therein), a closely related problem to \eqref{eq:ex_policy_eval} is 
	\begin{align}\label{eq:policy evaluation problem}
	\min_{\bx\in\mR^n_+} D_{\text{KL}}(\mathbb{E}[\br_\varphi],\mathbb{E}[A^+_\xi]\bx),
	\end{align}
	where $A^+_\xi \triangleq \max[(I-\gamma\hat{P})\Phi,\bo_{|\cY|\times n}]$ and $D_{KL}(\cdot,\cdot)$ is the Kullback-Liebler (KL) divergence~\cite{csiszar1991least} (see \eqref{eq:kl_dist} for the definition). This problem is known as the \emph{stochastic Poisson linear inverse problem} (see, e.g., \cite{bauschke2017descent} -- Section~5.1). The KL distance is a natural measure that corresponds to noise of Poisson type and does not admit globally Lipschitz continuous gradient~\cite{bauschke2017descent}.
\end{example}

\begin{example}[Risk-averse optimization]\label{ex:risk_averse}
	Consider the mean-variance risk-averse optimization problem~\cite{ruszczynski2009risk}
	\begin{align}\label{eq:mean_var}
	\min_{\bx\in\cX\subseteq\mR^n} -\mathbb{E}\left[r_\xi(\bx)\right] + \lambda \mathbb{V}\text{ar}[r_\xi(\bx)],
	\end{align}
	where $r_\xi(\bx):\mathbb{R}^n \rightarrow \mathbb{R}$ is a random function given by the decision variable $\bx$, $\mathbb{V}\text{ar}[\cdot]$ denotes the variance with respect to $\xi$, and $\lambda > 0$ is the \emph{risk-aversion} parameter. In the context of risk-averse optimal control of systems governed by partial differential equations (PDEs) under uncertainty, there is a need to solve \eqref{eq:mean_var} where $r_\xi(\bx)$ is the control objective which its evaluation requires the solution of a system of PDEs (see, e.g., \cite{borzi2009multigrid} and references therein). To render the computation of the control objective and its gradient tractable, it is common to employ quadratic approximations of the objective function in the context of risk-neutral optimal control~\cite{alexanderian2017mean}. However, the variance component of the risk-averse mean-variance formulation \eqref{eq:mean_var} results in \emph{polynomial objective function of degree 4} which does not have Lipschitz continuous gradient (see Section~\ref{sec:risk_averse_ex} for more details).
\end{example}

\begin{example}[Low-rank non-negative matrix factorization]
Let $D\in\mR^{p\times q}_{++}$ be an unknown non-negative matrix. We aim to find a low-rank approximate non-negative matrix factorization (NMF) of $D$ based on random samples $\hat{D}$ from a stochastic oracle such that $\mE[\hat{D}]=D$. This requires solving
\begin{align}\label{eq:low_rank}
	\min_{X\in\mR^{p\times k}_{+}, Y\in\mR^{q\times k}_{+}} \dist (\mE[\hat{D}],XY^\top),
	\end{align}
	where $k<\min\{p,q\}$ is a prespecified rank and $\dist (\cdot,\cdot)$ is some distance function~\cite{berry2007algorithms}. Low-rank approximate NMF has a range of applications in feature extraction, text mining, and spectral data analysis~\cite{lee1999learning,paatero1994positive}. Due to the non-negativeness of $D$ and $XY^\top$, there are information theoretic arguments on taking statistical divergences such as KL-divergence~\cite{lee1999learning} or $\varphi$-divergence~\cite{cichocki2006csiszar} for the distance function which do not have Lipschitz continuous gradient.
\end{example}

It is natural to use stochastic gradient descent (SGD) to solve \eqref{eq:2-level composition}. However, obtaining an \textit{unbiased} estimator of the gradient of the composition is challenging as
$$   \grad F(\bx) = \mathbb{E}_{\varphi,\xi}\left[  \grad g_\xi(\bx)\grad f_\varphi(\mathbb{E}_{\xi}[g_\xi(\bx)])  \right] \neq \mathbb{E}_{\varphi,\xi}\left[ \grad g_\xi(\bx)\grad f_\varphi(g_\xi(\bx))   \right],$$ i.e., the stochastic gradient of the outer function evaluated at a stochastic value of the inner function results in a biased estimator.

To overcome this challenge, different methods are proposed to approximate the inner expectation which then helps to prove convergence to an approximate stationary point. These methods are discussed in Section~\ref{sec:related_work}. However, similar to gradient-based methods for differentiable functions, theoretical guarantees for these methods heavily make use of the Lipschitz continuity of the gradient which narrows their applications. In the absence of Lipschitz continuity of the gradient, a sequence of recent works \cite{bauschke2017descent,lu2018relatively,bolte2018first} introduce the notion of \textit{relative smoothness} which provides a new non-quadratic upper bound to the objective function and provides a new descent lemma~\cite{bauschke2017descent}. The definition of relative smoothness is provided in Section~\ref{sec:prelem}. The new upper bound requires solving a Bregman gradient subproblem and is discussed in Section~\ref{sec:related_work}.

\subsection{Related Work}\label{sec:related_work}
This work is on the intersection of two research areas: 1. stochastic composition optimization, 2. optimization of functions without Lipschitz continuous gradient.
Below, we review the literature corresponding to these areas.

\textbf{Stochastic composition optimization.} 
The stochastic composition problem dates back to \cite{ermoliev1976stochastic} and regained attention due to the broad range of applications in machine learning. Large deviation bounds for the empirical optimal value were established in \cite{ermoliev2013sample}, and the central limit theorem for the composition problem was established in \cite{dentcheva2017statistical} which proved that the empirical optimal value of \eqref{eq:2-level composition} converges to the true optimal value. Recently, different works focus on developing efficient first-order algorithms to solve \eqref{eq:2-level composition}. \cite{wang2017stochastic} proposed a two-step algorithm to solve two-level composition problems. The first step approximates the inner function value by moving average
\vspace{-0.2cm}
\begin{align*}
	y^{k+1} = (1-\beta_k)y^k + \beta_k  {g_{\xi_k}}(\bx^k),
\end{align*}
with $\beta_k\in (0,1),$ and the second step performs the projected stochastic gradient update
\vspace{-0.2cm}
\begin{align*}
	\bx^{k+1} = \Pi_{\cX}\{\bx^k - \alpha_k \grad  {g_{\xi_k}}(\bx^k)\grad f_\varphi(y^{k+1})\},
\end{align*}
where $\Pi_\cX$ refers to the projection operator onto the set $\cX$.
Since the gradient estimator is biased, the algorithm requires the stepsize $\alpha_k$ to be small compared to $\beta_k$, i.e,  {it should satisfy $\lim_{k\rightarrow \infty} \alpha_k/\beta_k = 0$} which results in a \textit{two-timescale} algorithm. Almost sure convergence of the iterates of the algorithm to an optimal solution in the convex and to a stationary point in the nonconvex  {but unconstrained} setting are established \emph{under the smoothness assumption}. However, due to the two-timescale requirement of the algorithm, its sample complexity is inferior to SGD for single-level problems.  {The paper improves its sample complexity for the convex setting by a second algorithm.} Furthermore, the convergence rate of the two-timescale algorithm was improved in \cite{wang2016accelerating} where they further considered minimizing the sum of a two-level composition and a nonsmooth but convex regularization function. \cite{hu2020biased} also studied this problem under the same assumptions, but the sample complexity of their algorithm is still suboptimal compared to SGD.

Recently, \cite{ghadimi2020single} proposed a single-timescale algorithm (i.e., the one in which the rate of inner function and the outer function gradient updates are in the same order) to solve \eqref{eq:2-level composition}.  The paper also establishes almost sure convergence of the algorithm to a stationary point in the constrained setting. Motivated by gradient flow in continuous space, \cite{chen2020solving} proposes a single-timescale algorithm for the unconstrained problem. Both algorithms in \cite{ghadimi2020single} and \cite{chen2020solving} achieve the sample complexity of SGD. The difference between their algorithms is that \cite{ghadimi2020single} uses moving average on the inner function value, the $\{\bx^k\}$ sequence, and the gradient estimates, while
\cite{chen2020solving} improves the inner function value update of \cite{wang2017stochastic} as
\begin{align*}
	y^{k+1} = (1-\beta_k)(y_k + \grad  {g_{\xi_k}}(\bx^k)(\bx^k - \bx^{k-1})) + \beta_k  {g_{\xi_k}} (\bx^k),
\end{align*}
where $\grad g_\xi(\bx^k)(\bx^k - \bx^{k-1})$ is a correction term derived from an ODE analysis. The algorithm in \cite{chen2020solving} is also generalized for multilevel composition problems.

Besides \cite{chen2020solving}, there are other works that propose algorithms for  multilevel composition problems.  {The work of \cite{ghadimi2020single} is generalized in \cite{Balasubramanian2022} to solve multilevel composition problems with and without a need to mini-batch sampling in each iteration. \cite{Zhang21MultiLevel} also considers the multilevel composition problem  and proposes a variance reduction method to solve it. Total sample complexity of their method has polynomial dependence on the number of the levels.} \cite{yang2019multilevel} proposes a multi-timescale method and \cite{ruszczynski2021stochastic} generalizes the work of \cite{ghadimi2020single} to multilevel and, further, shows the almost sure convergence of their single-timescale algorithm for problems where the composition functions are Lipschitz continuous. The paper also provides the sample complexity of the algorithm when the composition functions are smooth.

The composition problem with specific structures are considered in a number of works. \cite{dai2017learning} considers the setting where the outer function $f$ is convex and the inner function $g$ is linear. 
 \cite{Blanchet17Unbiased} proposes an unbiased estimator of the gradient for convex problems in finite-sum form. For the same setting, \cite{lian2017finite,liu2018dualityfree,devraj2019stochastic,lin2020improved,xu2021katyusha,yu2017fast} use variance reduction techniques to accelerate the convergence. \cite{zhang2021stochastic} considers minimizing the sum of a composition and a nonsmooth convex function $R(\bx)$ and proposed a prox-linear algorithm with subproblem
\vspace{-0.2cm}
\begin{align*}
	\bx^{k+1} = \argmin _\bx f(\tilde{g}^k + \tilde{J}^k(\bx - \bx^k)) + R(\bx) + \frac{M}{2}\|\bx - \bx^k\|^2,
\end{align*}
where $\tilde{g}^k$ and $\tilde{J}^k$  are (mini- batch or variance-reduced) estimates of $g(\bx^k)$ and $\grad g(\bx^k)$, respectively. Their theory requires the outer function $f$ to be deterministic, convex, and possibly nonsmooth while $g$ be a smooth function.

Finally, we note that the composition problems are also related to biased gradient methods which focus on establishing the stationarity using biased gradients of the objective functions, not necessarily in the composition structure~\cite{hu2017analysis}.

\textbf{Mirror descent for relatively smooth function}.  The mirror descent, also known as Bregman proximal gradient method, was originally developed by Nemirovski and Yudin (see~\cite{nemirovskij1983problem}) to solve $\min_{\bx\in\cX}  f(\bx)$ in general Banach space where the primal and dual spaces are not isometric to each other -- see \cite{bubeck2014convex} for more information. Given a continuously differentiable and strictly convex  mirror function $h$, the mirror descent algorithm performs the update
$\bx^{k+1} = \grad h^{-1}(\grad h(\bx^k) - \alpha_k \grad f(\bx^k)$,
where $\grad h^{-1}$ is the inverse of the gradient map. Later \cite{beck2003mirror} proved that the mirror descent update is indeed the solution to the Bregman gradient subproblem
\vspace{-0.2cm}
\begin{align*}
	\bx^{k+1} = \argmin_{\bx \in \cX} \fprod{\grad f(\bx^k),\bx - \bx^k} + \frac{1}{\alpha_k}D_h(\bx,\bx^k),
\end{align*}
where $D_h(\bx,\bx^k) \triangleq h(\bx) - h(\bx^k) - \fprod{\grad h(\bx^k),\bx - \bx^k}$ is the Bregman distance.

From the optimization perspective, mirror descent allows obtaining convergence rates with significantly less dependence on the dimension of the ambient space $n$ in certain cases. For instance, consider minimizing a convex and $C_f$ Lipschitz continuous function $f(\bx)$ over the simplex $\cX\triangleq\{\bx \in \mathbb{R}^n: \sum_i^n x_i = 1, x_i\geq 0\}$. As discussed in \cite{beck2003mirror}, the iteration complexity of the subgradient projection method for this problem is
$\min_{1\leq s \leq k} f(\bx^s) - \min_{\bx\in \cX} f(\bx) \leq \cO(1)C_f\sqrt{n}/\sqrt{k}$, while the complexity of the mirror descent algorithm with $h(\bx)=2^{-1}\norm{\bx}_p^2$ and $p=1+(\log n)^{-1}$ is
$\cO(1)\sqrt{\log n}/\sqrt{k}$, which improves the subgradient method by $\sqrt{n/\log n}$. Note that the resulting mirror descent algorithm has closed-form solution for its subproblem. The mirror descent algorithm is used in (stochastic) convex \cite{beck2003mirror,ben2001ordered,auslender2009projected,nemirovski2009robust,duchi2010composite,juditsky2011first,asi2019modeling} and nonconvex \cite{ghadimi2016mini,wang2019spiderboost} settings. In either case, similar to other first-order methods, the Lipschitz continuity of the objective function or its gradient are crucial for the analysis.


In many problems, however, the objective function is not smooth, i.e., it does not have Lipschitz continuous gradient. Recently, a new \textit{relative smoothness} condition was introduce by \cite{bauschke2017descent} for convex function $f$, which generalizes the commonly used smoothness condition.
The function $f$ is smooth relative to $h$ if $Lh - f$ is a convex function -- see Definition~\ref{def:rel_smooth} for the formal definition. This condition provides an upper bound on the objective function through the Bregman distance with the generating function $h$ as
$f(\by) \leq f(\bx) + \fprod{\grad f(\bx),\by-\bx} + LD_h(\by,\bx)$,
which requires solving such problems using the mirror descent update discussed above. Latter, \cite{bolte2018first} generalized this smoothness condition for nonconvex functions and named it smooth adaptable (smad) property. \cite{lu2018relatively} studied the same property with a different requirement on the generating function $h$. This property unifies the mirror descent and the gradient descent as setting $h(\bx) = \frac{1}{2}\norm{\bx}^2$ results in $D_h(\by,\bx) = \frac{1}{2}\norm{\by-\bx}^2$. If so, the above inequality recovers the traditional smoothness inequality
$f(\by) \leq f(\bx) + \fprod{\grad f(\bx),\by-\bx} + \frac{L}{2}\norm{\by-\bx}^2$ 
and the corresponding descent lemma~\cite{teboulle2018simplified}.
Further studies focus on nonconvex \cite{bolte2018first,bauschke2019linear,mukkamala2020convex,ahookhosh2021bregman,li2019provable} and stochastic \cite{zhang2018convergence,dragomir2021fast,hanzely2021fastest,davis2018stochastic,hanzely2021accelerated} scenarios. The lower bound on the sample complexity under this new smoothness assumption was established in \cite{dragomir2021optimal} which indicates the need for additional assumptions for acceleration. Finally, we note that the notion of Lipschitz continuity was also generalized to relative continuity in \cite{davis2018stochastic,teboulle2018simplified}.

\subsection{Contributions}
This paper proposes three algorithms to solve the two-level stochastic composition problem \eqref{eq:2-level composition} and investigates their iteration/sample complexities in the absence of the smoothness of the inner, outer, or both composition functions. More specifically, the contributions of the paper are as follows.

\begin{enumerate}
	\item We consider three different scenarios of constrained two-level stochastic composition problems: Smooth of Relative-smooth(SoR), Relative-smooth of Smooth(RoS), and Relative-smooth of Relative-smooth(RoR). For each combination, we prove that the composition function is smooth relative to a (strictly/strongly) convex distance generating function.
	\item For the SoR composition, we propose a single-timescale mini-batch stochastic composition Bregman gradient method. This algorithm achieves an $\epsilon$-approximate stationary solution in $\cO(\epsilon^{-2})$  calls to the inner function value and gradient and outer function gradient stochastic oracles, which matches the sample complexity of SGD in the single-level problems. This result also matches the lower bound of the Bregman method with merely relative smoothness assumption \cite{dragomir2021optimal}.
	\item For the RoR (and RoS) composition,  we propose a mini-batch prox-linear stochastic Bregman gradient algorithm with  stochastic oracle complexity of $\cO(\epsilon^{-2})$ for the gradient of the inner and outer functions, and stochastic oracle complexity of $\cO(\epsilon^{-3})$ for the inner function value.
	\item For the RoS composition, we improve the RoR algorithm by a variance reduction technique in  {the RoS-VR algorithm}. The algorithm improves the stochastic oracle complexity of the inner function value to $\cO(\epsilon^{-5/2})$ and the inner function gradient to $\cO(\epsilon^{-3/2})$, while the complexity of the outer function gradient stochastic oracle is still $\cO(\epsilon^{-2})$.
\end{enumerate}

\vspace{-0.4cm}
\subsection*{Preliminaries}\label{sec:prelem}
Given a continuously differentiable function $f:\mR^n\rightarrow\mR^d$, $\grad f(\bx)\in\mR^{n\times d}$ denotes its Jacobian (equivalently its gradient when $d=1$). For a real-valued two-time continuously differentiable function $f$, $\grad^2 f$ denotes its Hessian. Given a closed convex set $\cX\subseteq\mR^n$, $\mathbf{int}\cX$ denotes its interior, $\cN_{\cX}(\bx)$ denotes the normal cone to set $\cX$ at $\bx$, and $\delta_{\cX}(\cdot)$ denotes the indicator function of $\cX$. Finally, given $\bx\in\mR^n$, $\norm{\bx}$ denotes its $l_2$ norm.

Next, we introduce the Bregman distance and the notion of relative smoothness which are central to the ideas discussed in this paper.

\begin{definition}[Bregman Distance]
	Let $h$ be a proper, differentiable and convex function on the $\mathbf{int}\dom h$. The Bregman distance of $\bx,\by\in \dom h$ generated by $h$ is defined as
	\vspace{-0.1cm}
	\begin{equation}
		D_h(\bx,\by) \triangleq  h(\bx)-h(\by)-\langle \grad h(\by),\bx-\by \rangle.
	\end{equation}
\end{definition}
Note that since $h$ is convex, $D_h(\bx,\by) \geq 0$ for any $\bx,\by \in \dom h$.  Furthermore, if $h$ is $\mu$-strongly convex, i.e., $h(\bx)-h(\by)-\langle \grad h(\by),\bx-\by \rangle \geq (\mu/2)\norm{\bx-\by}_2^2$, then $D_h(\bx,\by)\geq (\mu/2)\norm{\bx-\by}_2^2$ by the definition. Note that the Bregman distance is convex with respect to the first argument. For $h(\bx) = \frac{1}{2}\|\bx\|^2$, $D_h(\bx,\by) = \frac{1}{2}\|\bx - \by\|^2$, i.e., the Euclidean distance.

\begin{definition}[Relative Smoothness \cite{bauschke2017descent}]\label{def:rel_smooth}
	Assume $\dom h \subset \dom f$. The function $f$ is $L$-smooth relative to the convex function $h$ if
	\vspace{-0.2cm}
	\begin{equation} \label{eq:rel_smooth-univar}
		|f(\bx) - f(\by) - \langle \grad f(\by),\bx-\by \rangle | \leq LD_h(\bx,\by), \quad \forall  \bx,\by \in\dom h.
	\end{equation}
	If $f$ is a vector-valued function, then it is $L$-smooth relative to $h$ if
	\vspace{-0.2cm}
	\begin{equation}
		\norm{f(\bx) - f(\by) - \langle \grad f(\by),\bx-\by \rangle} \leq LD_h(\bx,\by), \quad \forall  \bx,\by \in\dom h.
	\end{equation}
\end{definition}
If $h(\bx) = \frac{1}{2}\|\bx\|^2$, the relative smoothness results into the famous upper/lower bound $|f(\bx) - f(\by) - \langle \grad f(\by),\bx-\by \rangle |\leq L\norm{\bx-\by}_2^2$  provided that $f$ is smooth, i.e., $\norm{\grad f(\bx) - \grad f(\by)}\leq L\norm{\bx-\by}$. The left-hand-side inequality from \eqref{eq:rel_smooth-univar}, i.e., $-LD_h(\bx,\by) \leq f(\bx)-f(\by) - \langle \grad f(\by),\bx-\by \rangle$ is also known as \textit{L-weak convexity of $f$ relative to $h$}~\cite{zhang2018convergence}, which generalizes the notion of \textit{weak convexity}, i.e., $-(L/2)\norm{\bx-\by}_2^2 \leq f(\bx)-f(\by) - \langle \grad f(\by),\bx-\by \rangle$.

Some examples of relatively smooth functions are provided below.
\begin{example}[Relatively smooth functions]\label{eg:relative-smooth}\
	\begin{itemize}
		\item Let $f(\bx)$ be twice differentiable and the operator norm of its Hessian with respect to the $l_2$-norm satisfies $\|\grad^2 f(\bx)\| \leq p_r(\|x\|)$, where $p_r(\cdot)$ is a univariate $r$-th degree polynomial that satisfies $p_r(\alpha)\leq L(1+\alpha^r)$ for $\alpha\geq 0$ and $L>0$. Then, $f(\bx)$ is $L$-smooth relative to $h(\bx) = \frac{1}{r+2} \|\bx\|^{r+2} + \frac{1}{2}\|\bx\|^2$~\cite{lu2018relatively}.
		\item In the Poisson linear inverse problem, given a matrix $A\in \mathbb{R}_+^{m\times n}$, one needs to recover the signal $\bx\in \mathbb{R}_{++}^n $ from noisy observations $\bb$  such that $A\bx \approx \bb$ by minimizing the Kullback-Leibler divergence
		\begin{align}\label{eq:kl_dist}
		D_{KL}(\bb,A\bx) \triangleq  \sum_{i=1}^n ( b_i \log (b_i/(A\bx)_i+(A\bx)_i - b_i ,
		\end{align}
		which is $L$-smooth relative to $h(\bx) = -\sum_{i=1}^n \log x_i$ for $L\geq \|\bb\|_1 $ -- see~\cite{bauschke2017descent}.
	\end{itemize}
\end{example}

Next, we present two lemmas which are key to the analysis of the Bregman gradient methods.
\begin{lemma}[Three-point identity, (Lemma 3.1 in \cite{chen1993convergence})]
	Suppose $h$ is a differentiable and convex function. For any three points $\by, \bz \in\mathbf{int} \dom h$ and $\bx \in \dom h$, the following  identity holds,
	\vspace{-0.2cm}
	$$
		D_h(\bx,\by) + D_h(\by,\bz) = D_h(\bx,\bz) + \fprod{\grad h(\bz) - \grad h(\by),\bx-\by}.
	$$
\end{lemma}
This identity is the generalization of Pythagorean identity in Euclidean geometry.
The next lemma follows from the three-point identity above.

\begin{lemma}[Three-point inequality] \label{lemma:three-point ineq}
	Let $\varphi$ be a proper, lower semicontinuous and convex function,  $h$ be a differentiable convex function. Given $\tau>0$ and $\bx\in \cX \subset \mathbf{int}\dom h$, let
	$$\bx^+ \in  \argmin_{\by \in \cX} \varphi(\by) + \frac{1}{\tau}D_h(\by,\bx), $$
	which is unique if $h$ is strongly convex. Then, $\bx^+$  satisfies
	$$ \tau (\varphi(\bx^+) - \varphi(\by)) \leq D_h(\by,\bx) - D_h(\bx^+,\bx) - D_h(\by,\bx^+), \quad\forall \by \in \cX.$$
\end{lemma}

\begin{remark}
	In particular,  if we set $\varphi(\by) = \fprod{\bw,\by - \bx}$, then
	$$\tau\fprod{\bw,\bx^+ - \by} \leq  D_h(\by,\bx) - D_h(\bx^+,\bx) - D_h(\by,\bx^+), \quad\forall \by \in \cX.$$
\end{remark}


\section{Smooth of Relative-smooth(SoR) Composition}
As we deal with the composition of two functions, under different smoothness assumptions, there are four possible composition scenarios. The smooth of  smooth composition is the simplest one which results in the composition function to be also smooth and it is well studied in the literature of stochastic composition problem listed above. In this section, we study the Smooth of Relative-smooth (SoR) composition. First, we give the definition of the stationarity measure under the Bregman distance.

\begin{lemma}[Stationarity Measure]\label{lemma:SoR-stationarity measure}
	Given a $\mu_h$-strongly convex function $h$, define
	\begin{align}\label{eq:opt_obj_1}
		\hat{\bx}^+ \triangleq \argmin_{\by\in \cX} \fprod{\grad F(\bx),\by-\bx} + \frac{1}{\tau}D_h(\by,\bx),
	\end{align}
	where $\tau>0$. Then $\hat{\bx}^+ = \bx$ if and only if $-\grad F(\bx)\in \cN_\cX(\bx)$.
\end{lemma}

{\it Proof}	By the optimality condition, we have
\begin{align}\label{eq:breg_subproblem}
	0 \in \grad F(\bx) + \frac{1}{\tau}(\grad h(\hat{\bx}^+) - \grad h(\bx)) + \partial \delta_\cX(\hat{\bx}^+).
\end{align}
Hence, if $\hat{\bx}^+ = \bx$, then $-\grad F(\bx)\in \cN_\cX(\bx)$.
Next, assuming $-\grad F(\bx)\in \cN_\cX(\bx)$, we have
$\fprod{\grad F(\bx),\by-\bx} \geq 0, \forall \by \in \cX.$
By the strong convexity of $h$, $D_h(\by,\bx)\geq  (\mu_h/2)\norm{\by-\bx}^2\geq 0$, $\forall \by \in \cX$, with the equality when $\by = \bx$.  Hence, the minimum of \eqref{eq:opt_obj_1} is unique with the solution $\hat{\bx}^+ = \bx$.
\qed

The above Lemma shows that $\dist(\hat{\bx}^+,\bx)$ is a suitable measure for the stationarity. Indeed, similar measures are commonly used in constrained optimization literature, e.g., in \cite{ghadimi2016mini}, the distance function $\dist(\hat{\bx}^+,\bx) \triangleq \frac{1}{2\tau^2}\norm{\hat{\bx}^+ - \bx}^2$ is used as the stationarity measure which is equal to $\norm{\grad F(\bx)}^2$ when $\cX = \mathbb{R}^n$.
 Following Lemma~\ref{lemma:SoR-stationarity measure}, we use
$\dist(\hat{\bx}^{k+1},\bx^k) = D_h(\hat{\bx}^{k+1},\bx^k)/\tau^2$ to measure stationarity in the SoR composition.

\vspace{-0.2cm}
\subsection{Smoothness of the Composition Function}
Below, we provide the main assumptions for the SoR composition. These assumptions are common in stochastic (composition) optimization with the main difference being generalization of the smoothness of the inner function to relative-smoothness.
\begin{assumption}[Smooth of Relative-smooth (SoR) composition]\label{assumption:SoR}
	The functions $f_\varphi$ and $g_\xi$ satisfy the following conditions:
	\begin{enumerate}[label=(\alph*)]
		\item The function $g_\xi$ is average $L_g$-smooth relative to 1-strongly convex function $h_g$, i.e., $\forall \bx_1,\bx_2 \in \dom g$,
		      $$\mathbb{E}_\xi[\norm{g_\xi(\bx_1) - g_\xi(\bx_2) - \fprod{\grad g_\xi(\bx_2),\bx_1 - \bx_2}}] \leq L_g D_{h_g}(\bx_1,\bx_2).$$
		\item The function $f_\varphi$ is average $L_f$-smooth, i.e., $\forall \bu_1,\bu_2 \in \dom f$,
		      \vspace{-0.2cm}
		      $$\mathbb{E}_\varphi[\norm{\grad f_\varphi(\bu_1) - \grad f_\varphi(\bu_2)}^2] \leq L_f^2 \norm{\bu_1 - \bu_2}^2.$$
		\item The stochastic gradients of $f_\varphi$ and $g_\xi$ are bounded in expectation, i.e., $\forall \bu \in \dom f$ and $\forall \bx \in \dom g$,
		      \begin{align*}
			      \mathbb{E}_\varphi[\norm{\grad f_\varphi(\bu)}^2] & \leq C_f^2,  \quad 			      \mathbb{E}_\xi[\norm{\grad g_\xi(\bx)}^2]       \leq C_g^2.
		      \end{align*}
	\end{enumerate}
\end{assumption}
\begin{remark}\label{rem:stoch_bound_lip}
	By Jensen's inequality, Assumption~\ref{assumption:SoR}- { $(a\sim c)$} also hold for $f$ and $g$ as well, i.e., $g$ is $L_g$-smooth relative to $h_g$, $f$ is $L_f$-smooth, and $f,g$ are $C_f,C_g$-Lipschitz continuous.
	Furthermore, note that bounded stochastic gradient implies bounded variance, i.e., $\forall \bu \in \dom f$, we have
	\begin{align*}
		\mathbb{E}_\varphi[\|\grad f_\varphi(\bu^k) - \grad f(\bu^k)\|^2 ] & = \mathbb{E}_\varphi[\|\grad f_\varphi(\bu^k)\|^2] - \|\grad f(\bu^k)\|^2 \\
		                                                                   & \leq \mathbb{E}_\varphi[\|\grad f_\varphi(\bu^k)\|^2] \leq C_f^2.
	\end{align*}
	Similarly, $\mathbb{E}_\xi[\|\grad g_\xi(\bx^k) - \grad g(\bx^k)\|^2 ] \leq C_g^2$,  $\forall \bx \in \dom g$.
\end{remark}

Under Assumption~\ref{assumption:SoR}, the SoR composition is proved to be smooth relative to the generating function $h$ defined in Lemma~\ref{lemma:SoR smooth}, with its proof  in  Appendix~\ref{pf:SoR smooth}.

\begin{lemma}\label{lemma:SoR smooth}
	Under Assumption \ref{assumption:SoR}, $F(\bx) = f(g(\bx))$ is 1-smooth relative to $C_g^2L_f$-strongly convex function $h(\bx) = \frac{C_g^2L_f }{2}\|\bx\|^2 + C_f L_g h_g(\bx)$.
\end{lemma}

\begin{remark}\label{remark:RoS to SoS}
	If $h_g(\bx) = \frac{1}{2}\|\bx\|^2$,  {i.e., $h_g$ is also smooth,} then under Assumption~\ref{assumption:SoR}, the composition function $F(\bx)$ is $(C_g^2L_f+C_fL_g)$-smooth  {which matches the results derived in \cite{ghadimi2020single}.}
\end{remark}

\subsection{Proposed Algorithm for the SoR Composition}
By relative smoothness of the composition function, it is natural to solve this problem using the Bregman gradient method~\cite{bolte2018first}. But, as discussed in the introduction, it is impossible to get an unbiased estimate of the gradient of the objective function. Given the current point $\bx^k$, one needs to estimate the inner function value and then obtain a stochastic gradient of the outer function at this estimate. This requires tracking three random sequences: the iterate $\{\bx^k\}$,  inner function value estimate $\{\bu^k\}$ at the current iterate, and the gradient estimate of the composition function $\{\bw^k\}$.  These random sequences are defined on the probability space $(\Omega,\cF,P)$, where $\Omega$ is the sample space, $\cF_k$ is the $\sigma$-algebra generated by the algorithm up to iteration $k$, i.e.,
$$\{\bx^0,\dots,\bx^{k+1}, \bu^0,\dots, \bu^{k}, \bw^0, \dots, \bw^{k}\},$$
and $P$ is the probability measure.

Given an iterate $\bx$, an estimate of the the inner function value $\bu$, and i.i.d. sample batches $\cB_{g}$, $\cB_{\grad g}$, and $\cB_{\grad f}$, the mini-batch estimate of the inner function and the mini-batch estimate of the gradients of the inner and outer functions are calculated as
\begin{subequations}\label{eq:vsw_update}
	\begin{align}
		\bu(\bx;\cB_{g})       & = \frac{1}{\abs{\cB_g}}\sum_{\xi \in \cB_g} g_\xi(\bx),        \label{eq:u_update}                                                       \\
		\bv(\bx;\cB_{\grad g}) & = \frac{1}{\abs{\cB_{\mathsmaller{\grad} g}}}\sum_{\xi \in \cB_{\mathsmaller{\grad} g}}\grad g_\xi(\bx),        \label{eq:v_update}      \\
		\bs(\bu;\cB_{\grad f}) & =  \frac{1}{\abs{\cB_{\mathsmaller{\grad} f}}} \sum_{\varphi \in \cB_{\mathsmaller{\grad} f}} \grad f_\varphi (\bu). \label{eq:s_update}
	\end{align}
\end{subequations}
Note that not all of the above estimators are used in all of the algorithms presented later.

We propose Algorithm~\ref{alg:SoR} to solve the SoR composition and its sample complexity is derived in the rest of this section. To do so, beside Assumption~\ref{assumption:SoR}, we require two extra assumptions discussed below. Assumption~\ref{assump:SoR stochastic oracle} {-$(a)$} is a standard lower boundedness of the objective function for nonconvex optimization and  {$(b)$} is the boundedness of the variance of the stochastic inner function value which is necessary to deal with stochastic estimators.
\begin{assumption}\label{assump:SoR stochastic oracle}
	The functions $f_\varphi$, $g_\xi$, and $F$ satisfy the following conditions:
	\begin{enumerate}[label=(\alph*)]
		\item The function $F(\bx)$ is lower bounded by $F^*$.
		\item The variance of $g_\xi$ is bounded, i.e.,
		      \vspace{-0.2cm}
		      $$\mathbb{E}_\xi[\norm{g_\xi(\bx) - g(\bx)}^2]\leq \sigma_g^2, \quad \forall \bx \in \dom g.$$
	\end{enumerate}
\end{assumption}

\begin{algorithm}
	\caption{SoR algorithm  {(SoR})}\label{alg:SoR}
	\begin{algorithmic}[1]
		\REQUIRE  $\bx^0\in\cX, \bu^0\in\mR^d, \{\tau_k\}_{k\in\mN_+},\{\beta_k\}_{k\in\mN_+}\subset(0,1), h_g$
        \STATE Sample $\cB_{\grad f}^{ {0}}$, $\cB_{\grad g}^{ {0}}$ 
        and update $\bv^0 = \bv(\bx^0;\cB^{ {0}}_{\grad g})$, $\bs^0 = \bs(\bu^0;\cB^{ {0}}_{\grad f})$ using \eqref{eq:v_update}-\eqref{eq:s_update}
		\STATE Update $\bw^0=\bv^0\bs^0$
		\FOR{$k=0,1,\dots,K-1$,}
		\STATE Given $ h(\bx) =\frac{C_g^2 L_f}{2}\|\bx\|^2+C_fL_gh_g(\bx) $, solve
		\begin{equation}\label{eq:SoR subproblem}
			\bx^{k+1} = \argmin_{\by\in \cX} \fprod{\bw^k,\by - \bx^k} + \frac{1}{\tau_k}D_h(\by,\bx^k)\end{equation}
		\STATE Sample $\cB^{ {k+1}}_g$ and update
		\begin{equation}\label{eq:u update}
			\bu^{k+1} =\frac{1}{\abs{\cB_g^{ {k+1}}}}\sum_{\xi \in \cB_g^{ {k+1}}}\left[ (1-\beta_k)(\bu^k + g_{\xi}(\bx^{k+1}) - g_{\xi}(\bx^k)) + \beta_k g_{\xi}(\bx^{k+1}) \right]
		\end{equation}
		\STATE Sample $\cB^{ {k+1}}_{\grad g}$, $\cB^{ {k+1}}_{\grad f}$, update $\bv^{k+1} = \bv(\bx^{k+1};\cB^{ {k+1}}_{\grad g})$, $\bs^{k+1} = \bs(\bu^{k+1};\cB^{ {k+1}}_{\grad f})$ by \eqref{eq:v_update}-\eqref{eq:s_update}
		\STATE Update $\bw^{k+1}=\bv^{k+1}\bs^{k+1}$
		\ENDFOR
	\end{algorithmic}
\end{algorithm}
To prove the complexity of this algorithm, we first need to establish error bounds for the stochastic estimates of the inner function value and gradient of the composition function. The error of the inner function value estimate by $\bu^k$ is bounded in Lemma \ref{lemma:SoR-error of u} which is established by \cite{chen2020solving}. For the completion of the paper, we also give the proof in Appendix \ref{pf:SoR-error of u}.

\begin{lemma}[Lemma 1 in \cite{chen2020solving}]\label{lemma:SoR-error of u}Under  {Assumptions} \ref{assumption:SoR} and \ref{assump:SoR stochastic oracle}, given $\beta_k\in(0,1)$, sequences  $\{\bx^k\}$ and $\{\bu^k\}$ generated by Algorithm~\ref{alg:SoR} satisfy
	\begin{align*}
		     & \mathbb{E}[\|g(\bx^{k+1}) - \bu^{k+1}\|^2 | \cF_{k}]                                                                                                \\
		\leq & (1-\beta_k)^2\|g(\bx^k) - \bu^k\|^2 + 4(1-\beta_k)^2C_g^2\|\bx^{k+1} - \bx^k\|^2 + 2\beta_k^2\sigma_g^2 /\abs{\cB^{ {k+1}}_g}.\numberthis \label{eq:u var}
	\end{align*}
\end{lemma}

Following the bound in Lemma~\ref{lemma:SoR-error of u}, the error of the composition gradient estimate $\bw^k$ is bounded as shown in Lemma~\ref{lemma:SoR-error of w}; the proof is provided in Appendix \ref{pf:SoR-error of w}.

\begin{lemma}\label{lemma:SoR-error of w}
	Under Assumption \ref{assumption:SoR}, the sequences $\{\bx^k\}$ and $\{\bw^k\}$ generated by Algorithm \ref{alg:SoR} satisfy
	\begin{align}\label{eq:gradient variance}
		\mathbb{E}[\|\grad F(\bx^k) - \bw^k\|^2]
		\leq {2C_f^2C_g^2}/{|\cB_{\grad f}^{ {k}}|} + {2C_f^2C_g^2}/{|\cB_{\grad g}^{ {k}}|}+ 2C_g^2L_f^2\mathbb{E}[\|g(\bx^k) - \bu^k\|^2 ].
	\end{align}
\end{lemma}

We investigate the sample complexity of the SoR algorithm by analyzing the merit function
\begin{equation}
	V(\bx^k,\bu^k) \triangleq  F(\bx^k)-F^*+\|g(\bx^{k})-\bu^{k}\|^2,
\end{equation}
which is denoted by $V^k$ in the rest of the paper. The per iteration merit function decrease is essential for the analysis and is presented in Lemma \ref{lemma:SoR-Vk decrease} with its proof in Appendix \ref{pf:SoR-Vk decrease}.

\begin{lemma}\label{lemma:SoR-Vk decrease}
	Let $\{\bx^k,\bu^k,\bw^k\}$ be the sequence generated by Algorithm~\ref{alg:SoR}. Furthermore, let Assumptions~\ref{assumption:SoR} and \ref{assump:SoR stochastic oracle} hold. Then
{\small
 \begin{align*}
		 & \mathbb{E}[V^{k+1}|\cF_k] \leq V^{k} + \frac{\tau_k}{2C_g^2L_f}\|\grad F(\bx^k) -\bw^k\|^2 -\left(\frac{1}{\tau_k}-2\right)D_h(\hat{\bx}^{k+1},\bx^k)                                   \\
		 & -\left(\frac{1}{\tau_k}-1 - \frac{8}{L_f} \right)D_h(\bx^{k+1},\bx^k) +((1-\beta_k)^2-1)\|g(\bx^k)-\bu^k\|^2 + {2\beta_k^2\sigma_g^2}/{|{\cB^{ {k+1}}_g}|}, \numberthis \label{eq:bound}
	\end{align*}
 }
	where $\hat{\bx}^{k+1}$ is defined as
	$$ \hat{\bx}^{k+1} \triangleq  \argmin_{\by \in \cX} \fprod{\grad F(\bx^k),\by - \bx^k} + \frac{2}{\tau_k}D_h(\by,\bx^k).$$
\end{lemma}
This Lemma establishes the connection of the merit function $V^k$, two estimates' errors $\norm{g(\bx^k)-\bu^k}^2$ and $\norm{\grad F(\bx^k) -\bw^k}^2$, and the stationarity measure $D_h(\hat{\bx}^{k+1},\bx^k)$ - see Lemma~\ref{lemma:SoR-stationarity measure}. Hence, by properly choosing the parameters, we can derive the sample complexity of this algorithm to obtain an $\epsilon$-stationary solution, which is presented in Theorem~\ref{thm:SoR-rate of convergence}, with the proof in Appendix \ref{pf:SoR-rate of convergence}.

\begin{theorem}[Sample complexity of the SoR algorithm]\label{thm:SoR-rate of convergence} Let $\{\bx^k\}$ be a sequence generated by Algorithm~\ref{alg:SoR} with $0<\tau_k<\min\{1/2,L_f/(L_f+8)\}$ and $\beta_k \in (0,1)$ satisfying $(1-\beta_k)^2 + \tau_kL_f -1 \leq 0$. Then, under Assumptions \ref{assumption:SoR} and \ref{assump:SoR stochastic oracle}, we have
	{\small
		\begin{align*}
		 & \mathbb{E}\left[ {D_h(\hat{\bx}^{R+1},\bx^R)}/{\tau_R^2} \right]                              \leq  \frac{V^0}{\sum_{j=0}^{K-1}(\tau_j - 2\tau_j^2)}+\frac{\sum_{k=0}^{K-1} \frac{\tau_k C_f^2}{L_f\abs{\cB^{ {k}}_{\grad f}}} + \frac{\tau_k C_f^2}{L_f\abs{\cB^{ {k}}_{\grad g}}} + \frac{2\beta_k^2\sigma_g^2}{\abs{\cB_g^{ {k+1}}}} }{\sum_{j=0}^{K-1}(\tau_j - 2\tau_j^2)} , 	\numberthis \label{eq:SoR thm}
	\end{align*}}
where $V^0 = F(\bx^0) -F^*+ \norm{g(\bx^0) - \bu^0}^2 $, and the expectation is taken with respect to all  random sequences generated by the algorithm and an independent random integer $R \in\{0,\dots,K-1\}$ with probability distribution
	{\small
		\begin{align*}
		P(R=k) = \frac{\tau_k - 2\tau_k^2}{\sum_{j=0}^{K-1}(\tau_j - 2\tau_j^2)}.
	\end{align*}}
\end{theorem}

\begin{remark}
	When $\tau_k \in (0,1/L_f)$, with $\beta_k \in (0,1)$,  the requirement $(1-\beta_k)^2 + \tau_kL_f -1 \leq0$ results in
	$1-\sqrt{1-\tau_kL_f}\leq \beta_k < 1$.
	Furthermore, note that  $\beta_k = \tau_kL_f$ satisfies the above bound  as $\tau_k < \frac{1}{L_f}$, so $\tau_k L_f <1$, and
	{\small
		\begin{align*}
		\beta_k - (1-\sqrt{1-\tau_k L_f}) & = (1+\sqrt{1-\tau_k L_f})(1-\sqrt{1-\tau_k L_f}) - (1-\sqrt{1-\tau_k L_f}) \\
		                                  & = (1-\sqrt{1-\tau_k L_f})\sqrt{1-\tau_k L_f} \ > 0.
	\end{align*}
	}

\end{remark}

\begin{remark}
	If we set $\abs{\cB^k_{\grad f}} = \abs{\cB^k_{\grad g}}\equiv 1$  {for all $k \in \{0, \dots, K-1\}$}, then
	{\small
	\begin{align*}
		\frac{1}{\sum_{j=0}^{K-1}(\tau_j - 2\tau_j^2)}\sum_{k=0}^{K-1}\left(\frac{\tau_k C_f^2}{L_f\abs{\cB^{ {k}}_{\grad f}}} + \frac{\tau_k C_f^2}{L_f\abs{\cB^{ {k}}_{\grad g}}} \right) = \frac{2C_f^2}{L_f}\frac{\sum_{k=0}^{K-1}\tau_k}{\sum_{j=0}^{K-1}(\tau_j - 2\tau_j^2)} \geq \frac{2C_f^2}{L_f},
	\end{align*}
	}%
	and hence, without the mini-batch, the right hand side of \eqref{eq:SoR thm} can not be smaller than an arbitrary $\epsilon>0$, and the  upper bound of \eqref{eq:SoR thm} cannot be made small enough. 
\end{remark}
Fixing the step and batch sizes for all iterations, the sample complexity of the algorithm is presented in a more readable form in Corollary \ref{cor:SoR-rate of convergence}. The proof follows from Theorem~\ref{thm:SoR-rate of convergence}.
\begin{corollary}\label{cor:SoR-rate of convergence}
	Given the assumptions of Theorem~\ref{thm:SoR-rate of convergence}, setting $\tau_k \equiv \tau <\min\{1/2,L_f/(L_f+8),1/L_f\}, \beta_k \equiv L_f\tau $, $\abs{\cB^k_{\grad f}} = \abs{\cB^k_{\grad g}} \equiv \left\lceil \frac{4C_f^2}{(1-2\tau)L_f \epsilon} \right\rceil$, and $ \abs{\cB^k_{g}} \equiv \left\lceil \frac{4\tau L_f^2\sigma_g^2}{(1-2\tau)\epsilon}\right\rceil$,
	we have
	\begin{align*}
		\frac{1}{K}\sum_{k=0}^{K-1}\mathbb{E}\left[ {D_h(\hat{\bx}^{k+1},\bx^k)}/{\tau^2} \right] & \leq  \frac{V^0}{K(\tau - 2\tau^2)} +  \epsilon.
	\end{align*}
\end{corollary}
Hence to achieve $\epsilon$-stationarity, we can set $K = \cO(\epsilon^{-1})$, which means that we need $\cO(\epsilon^{-2})$ calls to the $g_\xi$, $\grad g_\xi $ and $\grad f_\varphi$ stochastic oracles.

Besides the stationarity error, we can show that the errors of the inner function estimate $\bu^k$ and the gradient of the composition estimate $\bw^k$ can also be bounded with the same rate, as shown in Corollary~\ref{cor:SoR-rate of other estimators}. The proof is in Appendix \ref{pf:SoR-rate of other estimators}.
\begin{corollary}\label{cor:SoR-rate of other estimators}
	Under the setting of Corollary \ref{cor:SoR-rate of convergence}, we have
	{\small
	\begin{align*}
		\frac{1}{K}\sum_{k=0}^{K-1}\mathbb{E}[\|g(\bx^k) - \bu^k\|^2]\leq & \frac{\|g(\bx^0)-\bu^0\|^2}{KL_f^2\tau^2} + \frac{8V^0}{KL_f^3(\tau - \frac{L_f + 8}{L_f}\tau^2)}           \\
		                                                                  & + \left(\frac{8(1-2\tau)}{L_f^3(1 - \frac{L_f + 8}{L_f}\tau)} + \frac{1-2\tau}{2L_f^2\tau} \right)\epsilon,
	\end{align*}
	and
	\begin{align*}
		\frac{1}{K}\sum_{k=0}^{K-1}\mathbb{E}[\|\grad F(\bx^k) - \bw^k\|^2] & \leq \frac{2C_g^2\norm{g(\bx^0)-\bu^0}^2}{K\tau^2}  + \frac{16C_g^2V^0}{K(L_f\tau-(L_f+8)\tau^2)}                                  \\
		                                                                    & \quad + \left(\frac{16C_g^2(1-2\tau)}{L_f - (L_f+8)\tau} + \frac{C_g^2(1-2\tau)}{\tau} + \frac{\tau - 2\tau^2}{2} \right)\epsilon.
	\end{align*}
	}%
\end{corollary}
Hence, after $\cO(\epsilon ^{-2})$ calls to the stochastic oracles, we have $\frac{1}{K}\sum_{k=0}^{K-1}\mathbb{E}[\|\bu^k - \grad g (\bx^k)\|^2] \leq \cO(\epsilon)$ and $\frac{1}{K}\sum_{k=0}^{K-1}\mathbb{E}[\|\grad F(\bx^k) - \bw^k\|^2] \leq \cO(\epsilon)$.



\section{Relative-smooth of Relative-smooth (RoR) and \\ Relative-smooth of Smooth (RoS) Compositions}
This section considers the other two possible composition scenarios, namely relative-smooth of relative-smooth (RoR) and relative-smooth of smooth (RoS) compositions. First, Section~\ref{sec:ror_ros} establishes the relative smoothness of the composition function for the RoR and RoS settings - Lemma~\ref{lemma:RoR-relative smoothness}. Next, Section~\ref{sec:ror} provides an algorithm that solves both RoR and RoS settings and analyzes its sample complexity. Finally, Section~\ref{sec:ros} provides a variance-reduced algorithm with improved sample complexity for the RoS composition.

\subsection{Relative-smoothness of the Composition}\label{sec:ror_ros}
We first introduce our assumptions. While Assumption~\ref{assumption:RoR} is for more general RoR composition, we also included the Assumption~\ref{assumption:RoS} for the RoS setting which is specifically used in Section~\ref{sec:ros} for the analysis of the variance-reduced algorithm.

\begin{assumption}[Relative-smooth of Relative-smooth (RoR) composition] \label{assumption:RoR}
	The functions $f_\varphi, g_\xi$, and $h_f$ satisfy the following conditions:
	\begin{enumerate}[label=(\alph*)]
		\item The function $g_\xi$ is average $L_g$-smooth relative to 1-strongly convex function $h_g$.
		\item The function $f_\varphi$ is average $L_f$-smooth relative to convex function $h_f$.
		\item The stochastic gradients of $f_\varphi,g_\xi$ are bounded in expectation.
		\item The function $h_f$ is $C_{h_f}$-Lipschitz continuous, i.e.,
		      $\|\grad h_f(\bu)\| \leq C_{h_f}, \forall \bu \in \dom h$.
	\end{enumerate}
\end{assumption}

\begin{assumption}[Relative-smooth of Smooth (RoS) composition] \label{assumption:RoS}
	The functions $f_\varphi, g_\xi$, and $h_f$ satisfy the following conditions:
	\begin{enumerate}[label=(\alph*)]
		\item The function $g_\xi$ is average $L_g$-smooth.
		\item The function $f_\varphi$ is average $L_f$-smooth relative to convex function $h_f$.
		\item The stochastic gradients of $f_\varphi,g_\xi$ are bounded in expectation.
		\item The function $h_f$ is $C_{h_f}$-Lipschitz continuous.
	\end{enumerate}
\end{assumption}

The following proposition from \cite{zhang2018convergence}  {is used in Lemma~\ref{lemma:RoR-relative smoothness}} to establish the convexity of the distance generating function  {of} the composition $f(g(\bx))$ by showing its relative weak convexity  {(relative weak convexity is defined in the preliminary below Definition~\ref{def:rel_smooth}).} Note that Lemma~\ref{lemma:RoR-relative smoothness} is written for the more general RoR setting and it covers the RoS composition by setting $h_g(\bx) = \frac{1}{2} \norm{\bx}^2$. The proof is provided in Appendix \ref{pf:RoR-relative smoothness}.

\begin{proposition}[Proposition 2.2(c) in \cite{zhang2018convergence}] \label{prop:composition_rwc}
	Let $X$ be a nonempty closed convex set, $f: \mathbb{R}^{d} \rightarrow \mathbb{R}$ is closed convex and $C_f$-Lipschitz continuous, and $g: \mathbb{R}^{n} \rightarrow \mathbb{R}^{d}$ is $L_{g}$-smooth relative to $h$. Then the composition $f\circ g:\mathbb{R}^n\rightarrow \mathbb{R}$ is $C_fL_g$-weakly convex relative to $h$.
\end{proposition}


\begin{lemma} \label{lemma:RoR-relative smoothness}
	Under Assumption \ref{assumption:RoR}, $F(\bx) $ is 1-smooth relative to $$h(\bx) = (C_fL_g+C_{h_f}L_fL_g)h_g(\bx) + L_fh_f(g(\bx)),$$ which is shown to be convex. Furthermore, if $h_g(\bx)$ is 1-strongly convex, then $h(\bx)$ is $C_fL_g$-strongly convex.
\end{lemma}

\subsection{Proposed Algorithm for RoR (and RoS) Composition}\label{sec:ror}
In this section, we present an algorithm for the RoR composition which can also be used for the RoS composition, as a special case.

The distance generating function $h(\bx)$ in Lemma~\ref{lemma:RoR-relative smoothness} contains the composition $h(g(\bx))$. As $g$ has the expectation form, one can only have access to stochastic approximation of the distance generating function $h$ which are indeed biased. Note that this is not the case in the SoR setting, as its distance generating function is independent of the stochastic oracle - see Lemma \ref{lemma:SoR smooth}. Therefore, in the RoR (and RoS) setting, we cannot evaluate $D_h(\by,\bx^k)$, exactly, which is needed to solve \eqref{eq:SoR subproblem}.

Inspired by the idea of \cite{zhang2021stochastic}, we consider to approximate $D_h(\by,\bx^k)$ by linearizing the inner function $g$. This linearization is supported by the following Lemma \ref{lemma:RoR-linearization bound}, which provides a new upper bound for $F(\by)$ when replacing $g(\by)$ by its linear approximation $g(\bx) + \fprod{\grad g(\bx),\by-\bx}$ inside the composite distance generating function. The proof is available in Appendix \ref{pf:RoR-linearization bound}.

\begin{lemma}\label{lemma:RoR-linearization bound}
	Under Assumption \ref{assumption:RoR}, $\forall \bx,\by \in \cX$, we have
	\begin{equation}
		h_f(g(\by))  \leq h_f(g(\bx)+\fprod{\grad g(\bx),\by-\bx})+C_{h_f} L_g D_{h_g}(\by,\bx).
	\end{equation}
	Furthermore, we have
	\begin{align*}
		F(\by)  \leq & F(\bx)+\fprod{\grad F(\bx) ,\by-\bx}  +L_f D_{h_f}(g(\bx)                      \\
		             & +\grad g(\bx)^\intercal(\by-\bx),g(\bx))+ \lambda D_{h_g}(\by,\bx),\numberthis
	\end{align*}
	where $\lambda \triangleq  C_fL_g + 2C_{h_f}L_fL_g$.
\end{lemma}

Lemma \ref{lemma:RoR-linearization bound} provides an upper bound to $F$ that can be iteratively minimized to obtain a stationary solution. If we have access to $\grad F(\bx^k),\grad g(\bx^k)$, and $g(\bx^k)$ or their approximations, then we can solve
\begin{align*}
	\tilde{\bx}^{k+1} =  \argmin_{\by\in\cX} & \langle \grad F(\bx^k) ,\by-\bx^k \rangle  +\frac{L_f}{\tau_k} D_{h_f}(g(\bx^k)+\grad g(\bx^k)^\intercal(\by-\bx^k),g(\bx^k))\\
	                                          & \quad + \frac{\lambda}{\tau_k} D_{h_g}(\by,\bx^k). \numberthis \label{eq:RoR-deterministic subproblem }
\end{align*}

The objective function of this subproblem is convex and always upper bounds $F(\tilde{\bx}^{k+1})$ if $\tau_k \leq 1$. However, given that functions $f$ and $g$ involve expectations, it is impossible to evaluate $\grad F(\bx^k), g(\bx^k)$, and $\grad g(\bx^k)$ in \eqref{eq:RoR-deterministic subproblem } exactly. Instead, Algorithm~\ref{alg:RoR-mini batch} proposes mini-batch estimates of these terms and then solves \eqref{eq:RoS/R subproblem} iteratively, to obtain a stationary solution. The $\cF_k$ for this Algorithm is generated from $$\{\bx^0,\cdots,\bx^k; \bu^0,\cdots,\bu^{k-1};\bv^0,\cdots,\bv^{k-1};\bs^0,\cdots,\bs^{k-1}\}. $$
\begin{algorithm}
	\caption{RoR mini-batch algorithm (RoR)}\label{alg:RoR-mini batch}
	\begin{algorithmic}[1]
		\REQUIRE  $\bx^0\in\cX, \tau_k <\min\left \{1, \frac{C_fL_g + 2C_{h_f}L_fL_g}{C_fL_g+2} \right\}, \lambda \triangleq  C_fL_g + 2C_{h_f}L_fL_g, h_f, h_g$
		\FOR{$k= 0, 1,\dots,K-1$}
		\STATE Sample $\cB_g^k$, $\cB_{\grad g}^k$, $\cB_{\grad f}^k$  and update
		$\bu^{k} = \bu(\bx^k;\cB_{g}^k)$, $\bv^{k} = \bv(\bx^{k};\cB^k_{\grad g})$, and \\ $\bs^{k} = \bs(\bu^{k};\cB^k_{\grad f})$ using \eqref{eq:u_update}-\eqref{eq:s_update}

		\STATE Update  $\bw^{k}=\bv^k\bs^k$

		\STATE Given $h_f$ and $h_g$ (with $h_g(\bx) = \frac{1}{2}\|\bx\|^2$ in the RoS case), solve
		\begin{align*}
			\bx^{k+1} =\argmin_{\by \in \cX} & \fprod{\bw^k,\by-\bx^k} + \frac{L_f}{\tau_k}D_{h_f}(\bu^k+ (\bv^k)^\intercal(\by - \bx^k),\bu^k) + \frac{\lambda}{\tau_k}D_{h_g}(\by,\bx^k) \numberthis \label{eq:RoS/R subproblem}
		\end{align*}

		\ENDFOR
	\end{algorithmic}
\end{algorithm}

 Lemma~\ref{lemma:RoR-stationarity measure} below motivates the stationarity measure used for the RoR and RoS compositions.
\begin{lemma}\label{lemma:RoR-stationarity measure}
	$\forall \alpha, \beta, \tau >0$, define
	{\small
		\begin{align}\label{eq:opt_obj_2}
		\tilde{\bx}^{+} \triangleq \argmin_{y\in\cX} \fprod{\grad F(\bx),\by-\bx} + \frac{\alpha}{\tau}D_{h_f}(g(\bx) + \grad g(\bx)^\intercal(\by-\bx),g(\bx)) + \frac{\beta}{\tau}D_{h_g}(\by,\bx).
	\end{align}
	}
	Then $\tilde{\bx}^{+} =\bx $ if and only if $-\grad F(\bx) \in \cN_{\cX}(\bx)$.
\end{lemma}
{\it Proof}
By the optimality condition, we have
\begin{align*}
	0 & \in \grad F(\bx) + \frac{\alpha}{\tau}\left(\grad g(\bx) \grad h_f(g(\bx) + \grad g(\bx)^\intercal (\tilde{\bx}^{+}-\bx)) - \grad g(\bx) \grad h_f(g(\bx)\right) \\
	  & \quad + \frac{\beta}{\tau}(\grad h_g(\tilde{\bx}^{+}) - \grad h_g(\bx)) + \cN_{\cX}(\tilde{\bx}^{+}).
\end{align*}
Hence, if $\tilde{\bx}^{+} = \bx$, we have $ -\grad F(\bx) \in  \cN_{\cX}(\bx)$.
For the other direction, from $ -\grad F(\bx) \in  \cN_{\cX}(\bx)$, we have  $\fprod{\grad F(\bx),\by-\bx}\geq 0$, $\forall \by \in \cX$. Also, by strong convexity of $h_g$, $D_{h_g}(\by,\bx) \geq 0$ where the equality holds only with $\bx = \by$. Given the fact $D_{h_f}(g(\bx) + \grad g(\bx)^\intercal (\by-\bx)$ and $g(\bx))\geq 0$, $\tilde{\bx}^{+} = \bx$ is the unique optimal solution to \eqref{eq:opt_obj_2}.
\qed

Following Lemma~\ref{lemma:RoR-stationarity measure}, we use
$\dist(\tilde{\bx}^{k+1},\bx^k)$ to measure stationarity in the RoR and RoS composition scenarios where $\tilde{\bx}^{k+1}$ is generated by the deterministic subproblem \eqref{eq:RoR-deterministic subproblem }.


The rest of this section provides the sample complexity of the RoR algorithm. First, Lemma~\ref{lemma:RoR-decrease lemma}  upper bounds the objective function at iteration $k+1$ with that of iteration $k$ and some estimation errors. The result is then used in Lemma~\ref{lemma:RoR sequence bounded} to upper bound the running average of the distance between two consecutive iterates.

\begin{lemma}\label{lemma:RoR-decrease lemma}
	Let $\{\bx^{k}\}$ be the sequence generated by Algorithm \ref{alg:RoR-mini batch}, under Assumption \ref{assumption:RoR}, we have
	\begin{align*}
		f(g(\bx^{k+1})) & \leq f(g(\bx^k)) + 2C_f\|g(\bx^k)-\bu^k\| + \frac{C_f^2}{2}\|\grad g(\bx^k) - \bv^k\|^2                                                                                       \\
		                & \quad + \frac{1}{2}\|\bv^k\|^2 \|\grad f(\bu^k)-\bs^k\|^2         - \left( \frac{\lambda}{\tau_k} - C_fL_g - 2 \right)D_{h_g}(\bx^{k+1},\bx^k), \numberthis \label{eq:RoR ub}
	\end{align*}
	where $\lambda = C_fL_g + 2C_{h_f}L_fL_g$.
\end{lemma}
The proof of the lemma is provided in Appendix \ref{pf:RoR-decrease lemma}. Fixing the step sizes and batch sizes for all iterations and using Lemma~\ref{lemma:RoR-decrease lemma}, the sample complexity of obtaining two consecutive iterates with distance less than $\epsilon$ on average (i.e., running average) is provided in Lemma~\ref{lemma:RoR sequence bounded}, with the proof in Appendix \ref{pf:RoR sequence bounded}.

\begin{lemma}\label{lemma:RoR sequence bounded}
	Let Assumptions \ref{assump:SoR stochastic oracle} and \ref{assumption:RoR} (or \ref{assumption:RoS})  hold. Then, by setting $\tau_k \equiv \tau < \min\left \{1, \frac{C_fL_g + 2C_{h_f}L_fL_g}{C_fL_g+2} \right\}$, $\abs{\cB_g^k} \equiv \left\lceil \frac{ 16C_f^2\sigma_g^2}{M_1^2\epsilon^2} \right\rceil$, and $\abs{\cB_{\grad f}^k} = \abs{\cB_{\grad g}^k} \equiv \left\lceil \frac{2C_f^2C_g^2}{M_1\epsilon} \right\rceil$,	we have
		\begin{align*}
		\frac{1}{K}\sum_{k=0}^{K-1} \mathbb{E} \left[{D_{h_g}(\bx^{k+1},\bx^k)}/{\tau^2}\right] \leq \frac{f(g(\bx^0)) - F^*}{M_1K} + \epsilon,
	\end{align*}
	where $M_1\triangleq (C_fL_g+2C_{h_f}L_fL_g)\tau - (C_fL_g+2)\tau^2 $.
\end{lemma}
Note that while $h_f$ is not smooth, we assume that it is smooth on \textit{any bounded subset} of $\mR^d$.

\begin{assumption}\label{assumption:hf-smooth}
	The functions $h_f$ and $f$ are twice differentiable. Furthermore,  $\grad h_f$ is $L_{h_f}$-Lipschitz continuous on any bounded subset of $\mathbb{R}^d$.
\end{assumption}

\begin{remark}
	The above assumption is also used in \cite{bolte2018first} to derive the convergence of their algorithm. Notice that following the setting of Lemma \ref{lemma:RoR sequence bounded}, for any  $\epsilon>0$ and finite $K$, the Bregman distance between two consecutive iterates of Algorithm \ref{alg:RoR-mini batch} is bounded in expectation; hence, the sequence $\{\bx^k\}_{k=0}^K$ lies in a bounded set with probability 1.
	Furthermore, since $g$ is Lipschitz continuous, $\{g(\bx^k)\}$ also lies in a bounded set with probability 1. With bounded $\{g(\bx^k)\}$ and  the assumption that the variance of $g_\xi$ is bounded, $\{\bu^k\}$ is also guaranteed to be bounded. Hence,  {the argument $\bu$ in $\grad h_f(\bu)$ in Assumption~\ref{assumption:hf-smooth} belongs to a bounded set} for any sequence generated by Algorithm \ref{alg:RoR-mini batch} with probability 1.
\end{remark}

\begin{lemma}\label{lemma:f-smooth}
	Under Assumptions~\ref{assumption:RoR} (or \ref{assumption:RoS}) and \ref{assumption:hf-smooth}, $f$ is $L_f L_{h_f}$-smooth on any bounded subset of $\mathbb{R}^d$.
\end{lemma}

The proof of Lemma \ref{lemma:f-smooth} is provided in Appendix \ref{pf:f-smooth}. Lemmas~\ref{lemma:RoR-relationship of x_tilde and x} and \ref{lemma:RoR wk bound} are needed to establish the upper bound on $\norm{\tilde{\bx}^{k+1}-\bx^k}^2$ in Theorem~\ref{thm:RoR-sample complexity}. The proof of Lemma~\ref{lemma:RoR-relationship of x_tilde and x} is provided in Appendix \ref{pf:RoR-relationship of x_tilde and x}, while the proof of Lemma~\ref{lemma:RoR wk bound} is similar to that of Lemma~\ref{lemma:SoR-error of w} and is omitted.

\begin{lemma}\label{lemma:RoR-relationship of x_tilde and x}
	Let $\tilde{\bx}^{k+1}$ be defined as in \eqref{eq:RoR-deterministic subproblem } with $\lambda = C_fL_g + 2C_{h_f}L_fL_g$. Under Assumptions \ref{assump:SoR stochastic oracle}, \ref{assumption:RoR} (or \ref{assumption:RoS}) and \ref{assumption:hf-smooth}, we have
 {\small
	\begin{align*}
		\|&\tilde{\bx}^{k+1} - \bx^k\|^2  \leq  \ \frac{2\tau_k^2}{C_f^2L_g^2}\|\bw^k - \grad F(\bx^k)\|^2  + \frac{8C_{h_f}L_f}{C_f L_g}\|\bu^k - g(\bx^k)\|                   \\
		                                    & + \frac{2C_g^2L_fL_{h_f}^2}{C_{h_f}C_fL_g^2}\|\bu^k - g(\bx^k)\|^2 + \frac{6C_{h_f}L_f}{C_fL_g^2}\|\bv^k - \grad g(\bx^k)\|^2 
		                                     + \frac{4(C_f + 3C_{h_f}L_f)}{C_f}D_{h_g}(\bx^{k+1},\bx^k). \numberthis \label{eq:x_tilde-bound}
	\end{align*}
 }
\end{lemma}

\begin{lemma}\label{lemma:RoR wk bound}
	Under Assumptions  \ref{assumption:RoR} (or \ref{assumption:RoS}) and \ref{assumption:hf-smooth}, the sequences $\{\bw^k\}$ and $\{\bx^k\}$ satisfy
	\begin{align*}
		\mathbb{E}[\|\grad F(\bx^k) - \bw^k\|^2|\cF^k ] \leq 2C_g^2L_f^2L_{h_f}^2\sigma_g^2/\abs{\cB_g^k} + 2C_f^2C_g^2/\abs{\cB^k_{\grad f}} + 2C_f^2C_g^2/\abs{\cB^k_{\grad g}}.
	\end{align*}
\end{lemma}

Finally, we can establish the sample complexity of  {the RoR algorithm}~\ref{alg:RoR-mini batch} to find an approximate stationary solution, which is presented in Theorem \ref{thm:RoR-sample complexity}, with the proof in Appendix \ref{pf:RoR-sample complexity}.
\begin{theorem}[Sample complexity of the RoR algorithm]\label{thm:RoR-sample complexity}
	Let $\{\bx^k\}$ be a sequence generated by Algorithm~\ref{alg:RoR-mini batch} with $\tau_k \equiv \tau < \min\left \{1, \frac{C_fL_g + 2C_{h_f}L_gL_f}{C_fL_g+2} \right\}$, $\abs{\cB_g^k} \equiv \left\lceil \frac{ 16C_f^2\sigma_g^2}{M_1^2\epsilon^2} \right\rceil$, and $\abs{\cB_{\grad f}^k} = \abs{\cB_{\grad g}^k} \equiv \left\lceil \frac{2C_f^2C_g^2}{M_1\epsilon} \right\rceil$,	where $M_1 =(C_fL_g+2C_{h_f}L_gL_f)\tau - (C_fL_g+2)\tau^2$. Then, under Assumptions \ref{assump:SoR stochastic oracle}, \ref{assumption:RoR} (or \ref{assumption:RoS}) and  \ref{assumption:hf-smooth}, we have
	{\small
	\begin{align*}
		 & \frac{1}{K}\sum_{k=0}^{K-1}\mathbb{E}\left[ {\norm{\tilde{\bx}^{k+1}-\bx^k}^2}/{\tau^2}\right]                                                                                                          \leq  \left(\frac{4C_f + 6C_{h_f}L_f }{M_1C_f}\right)\frac{f(g(\bx^0)) - F^*}{K} \\
		 & \qquad + \left(\frac{C_g^2L_f^2L_{h_f}^2M_1^2}{4C_f^4L_g^2} + \frac{C_g^2L_fL_{h_f}^2M_1^2}{8C_{h_f}C_f^3L_g^2\tau^2} \right)\epsilon^2                                                                                                                                                      \\
		 & \qquad + \left(\frac{4M_1}{C_f^2L_g^2} + \frac{2C_{h_f}L_fM_1}{C_f^2L_g\tau^2} + \frac{3C_{h_f}L_fM_1}{C_f^3L_g^2\tau^2} + 4 + \frac{12C_{h_f}L_f}{C_f}  \right)\epsilon.
	\end{align*}
	}
\end{theorem}
Hence, to obtain an $\epsilon$-stationarity solution, by setting $K = \cO(\epsilon^{-1})$, the RoR algorithm~\ref{alg:RoR-mini batch} requires  $\cO(\epsilon^{-3})$ calls to $g_\xi$ and $\cO(\epsilon^{-2})$ calls to $\grad f_\varphi$ and $\grad g_\xi$ stochastic oracles.

\subsection{Variance-reduced Algorithm for the RoS Composition}\label{sec:ros}
When the inner function is Lipschitz smooth, i.e., Assumption \ref{assumption:RoS}, the RoR algorithm can be improved by variance reduction to estimate the inner function and its gradient. We use the stochastic path integrated differential estimator to estimate $g$ and $\grad g$~\cite{fang2018spider,nguyen2017sarah}. These updates are incorporated in Algorithm~\ref{alg:RoS-SPIDER} - see \eqref{eq:u_update_spider} and \eqref{eq:v_update_spider}. The corresponding $\cF_j^k$ is now generated from
\begin{align*}
	\{\bx_0^0,\cdots,\bx^0_{J-1},\bx^1_0,\cdots,\bx^k_{j};\bu^0_0,\cdots,\bu^0_{J-2},\bu_0^1,\cdots,\bu_{j-1}^k; \\
	\bv^0_0,\cdots,\bv^0_{J-2},\bv_0^1,\cdots,\bv_{j-1}^k; \bs^0_0,\cdots,\bs^0_{J-2},\bs_0^1,\cdots,
	\bs_{j-1}^k\}.
\end{align*}
The error bounds for the new estimators are provided in the following lemma.
\begin{lemma}(Lemma 1 in \cite{fang2018spider})\label{lemma:RoS error}
	Let Assumptions \ref{assump:SoR stochastic oracle} and \ref{assumption:RoS} hold, and $\bu^k_j$ and $\bv^k_j$ be updated as in Algorithm \ref{alg:RoS-SPIDER}. Then, $\{\bu_j^k\}, \{\bv_j^k\}$, and $\{\bx^k_j\}$ satisfy
	\begin{align}
		\mathbb{E}[ \|\bu_j^k - g(\bx_j^k)\|^2]       & \leq \frac{\sigma_g^2}{\abs{\cB_g^k}} + \sum_{r=0}^{j-1} \frac{C_g^2}{|{\cS^{k,r+1}_g}|} \mathbb{E}[\|\bx_{r+1}^k - \bx_{r}^k\|^2 ],            \\
		\mathbb{E}[ \|\bv_j^k - \grad g(\bx_j^k)\|^2] & \leq \frac{C_g^2}{|{\cB_{\grad g}^k}|} + \sum_{r=0}^{j-1} \frac{L_g^2}{|{\cS^{k,r+1}_{\grad g}}|} \mathbb{E}[\|\bx_{r+1}^k - \bx_{r}^k\|^2 ].
	\end{align}
\end{lemma}

\begin{algorithm}[h]
	\small
	\caption{Variance-reduced RoS algorithm  {(RoS-VR)}}\label{alg:RoS-SPIDER}
	\begin{algorithmic}[1]
		\REQUIRE  $\bx^0_0\in\cX, \tau_k < 1, \lambda \triangleq  C_fL_g + 2C_{h_f}L_fL_g, h_f$
		\FOR{$k= 0, 1,\dots,K-1$}
		\FOR{$j = 0,1, \dots, J-1$}
		\IF{$j == 0$}
		\STATE Sample $\cB^k_g$ and use \eqref{eq:u_update} to update
		$\bu^{k}_0 = \bu(\bx^k_j;\cB_{g}^k)$

		\STATE Sample $\cB^k_{\grad g}$ and use \eqref{eq:v_update} to update
		$\bv^k_0=\bv(\bx^k_j;\cB_{\grad g}^k)$

		\ELSE
		\STATE Sample $\cS^{k,j}_g$ and calculate
		\begin{align}
			\bu^k_j = \bu^k_{j-1} + \frac{1}{|{\cS^{k,j}_g}|} \sum_{\xi \in \cS^{k,j}_g} \left(g_\xi(\bx^k_j) - g_\xi(\bx_{j-1}^k) \right) \label{eq:u_update_spider}
		\end{align}
		\STATE Sample $\cS^{k,j}_{\grad g}$ and calculate
		\begin{align}
			\bv^k_j = \bv^k_{j-1} + \frac{1}{|{\cS^{k,j}_{\grad g}}|} \sum_{\xi \in \cS^{k,j}_{\grad g}} \left(\grad g_\xi(\bx^k_j) - \grad g_\xi(\bx^k_{j-1}) \right) \label{eq:v_update_spider}
		\end{align}
		\ENDIF
		\STATE Sample $\cB^{k,j}_{\grad f}$ and used \eqref{eq:s_update} to update
		$\bs^k_j=\bs(\bu^k_j;\cB^{k,j}_{\grad f})$
		\STATE Update $\bw^k_j = \bv^k_j\bs^k_j$
		\STATE Solve
		{\scriptsize
		\begin{align*}
			\bx^{k}_{j+1} =\argmin_{\by \in \cX} & \fprod{\bw^k_j,\by-\bx^k_j} + \frac{L_f}{\tau_{k,j}}D_{h_f}(\bu^k_j+ (\bv^k_j)^\intercal(\by - \bx^k_j),\bu^k_j) + \frac{\lambda}{2\tau_{k,j}}\norm{\by-\bx^k_j}^2 \numberthis \label{eq:RoS subproblem}
		\end{align*}
		}%
		\ENDFOR
		\STATE Set $\bx_0^{k+1} = \bx_J^k$
		\ENDFOR
	\end{algorithmic}
	\normalsize
\end{algorithm}

Furthermore, as $\forall \delta>0$, we have
\begin{align*}
	\mathbb{E}[\|\bu_j^k - g(\bx_j^k)\| ] & \leq \sqrt{\frac{\sigma_g^2}{|{\cB_g^k}|}} + \sqrt{\sum_{r=0}^{j-1} \frac{C_g^2}{|{\cS_g^{k,r+1}}|} \mathbb{E}[\|\bx_{r+1}^k - \bx_{r}^k\|^2 ]}                                                                   \\
	                                      & \leq \frac{\sigma_g}{\sqrt{|{\cB_g^k}|}} + \frac{\delta}{2} + \frac{1}{2\delta} \sum_{r=0}^{j-1} \frac{C_g^2}{|{\cS^{k,r+1}_g}|} \mathbb{E}[\|\bx_{r+1}^k - \bx_{r}^k\|^2 ], \numberthis \label{eq:u-g_UB_spider}
\end{align*}
where the first inequality follows from Jensen's inequality and the fact that $\sqrt{a+b}\leq\sqrt{a}+\sqrt{b}$, and the second inequality is due to Young's inequality. These inequalities are used in the following lemma to provide a new error bound for the new estimator of the gradient of the composition.
\begin{lemma}\label{lemma:RoS-w error}
	Under Assumptions \ref{assump:SoR stochastic oracle} and \ref{assumption:RoS}, the sequences $\{\bw^k_j\}$ and $\{\bx^k_j\}$ generated by Algorithm \ref{alg:RoS-SPIDER} satisfy
	\begin{align*}
		\mathbb{E}[\|\grad F(\bx^k_j) - \bw^k_j\|^2 ] \leq & 2C_g^2L_f^2L_{h_f}^2\mathbb{E}[\|g(\bx^k_j) - \bu^k_j\|^2 ] + 2C_f^2\mathbb{E}[\|\grad g(\bx^k_j) - \bv^k_j\|^2 ] \\
		                                                   & + 2C_f^2C_g^2/|{\cB_{\grad f}^{k,j}}|.
	\end{align*}
\end{lemma}
The proof is similar to that of Lemma \ref{lemma:SoR-error of w} with adjustment of the indices and is omitted here.
We can now derive the sample complexity of Algorithm \ref{alg:RoS-SPIDER}. First, we will show that the average Bregman distance of two consecutive points are bounded (Lemma~\ref{lemma:RoS-sequence bound}, with the proof in Appendix \ref{pf:RoS-sequence bound}); hence, Assumption \ref{assumption:hf-smooth} is applicable. Next, under smoothness of $h_f$ on bounded subsets of $\mR^d$, we derive the bound on $D_h(\tilde{\bx}^k_{j+1},\bx^k_j)$ and the final sample complexity.

\begin{lemma}\label{lemma:RoS-sequence bound}
	Under Assumptions \ref{assump:SoR stochastic oracle} and \ref{assumption:RoS}, setting $\tau_k \equiv \tau < \frac{C_fL_g}{C_fL_g+2}$,  $\abs{\cS^{k,j}_g} \equiv \left \lceil\frac{6C_f^2C_g^2\tau}{C_{h_f}L_fL_g}\frac{J}{M_2\epsilon}\right\rceil$, $\abs{\cS^{k,j}_{\grad g}} \equiv \left\lceil \frac{3C_f^2L_g\tau J}{C_{h_f}L_f} \right \rceil$, $\abs{\cB^k_g}\equiv \left\lceil \frac{36C_f^2\sigma_g^2}{M_2^2\epsilon^2} \right \rceil$, and  $\abs{\cB^k_{\grad g}}\equiv \abs{\cB^{k,j}_{\grad f}}\equiv\left\lceil \frac{15C_f^2C_g^2}{2M_2\epsilon} \right \rceil$, where $M_2 \triangleq \frac{C_fL_g\tau}{2} - \left(\frac{C_fL_g}{2}+1\right)\tau^2$, the sequence $\{\bx^k_j\}$ 
 satisfies
	\begin{align*}
		\frac{1}{KJ} \sum_{k=0}^{K-1}\sum_{j=0}^{J-1}\mathbb{E}\left[{\|\bx^k_{j+1}-\bx^k_j\|^2}/{\tau^2} \right]\leq \frac{f(g(\bx^0_0)) -F^*}{M_2KJ} +  \epsilon.
	\end{align*}
\end{lemma}

\begin{theorem}[{Sample complexity of the  {RoS-VR} algorithm}]\label{thm:RoS-sample complexity}
	Under Assumptions \ref{assump:SoR stochastic oracle}, \ref{assumption:RoS} and \ref{assumption:hf-smooth}, following the setting of Lemma \ref{lemma:RoS-sequence bound}, the sequence $\{\bx^k\}$ generated by Algorithm \ref{alg:RoS-SPIDER} satisfies
	{\small
	\begin{align*}
		 & \quad\frac{1}{KJ}\sum_{k=0}^{K-1}\sum_{j=0}^{J-1} \mathbb{E}\left[\frac{\|\tilde{\bx}^k_{j+1}-\bx^k_j\|^2}{\tau^2} \right]                                                                                    \\
		 & \leq \left(\frac{2A_0}{C_fL_gM_1} + \frac{A_1C_{h_f}L_f\tau}{2C_f^2M_1^2} + \frac{2A_3C_{h_f}L_f\tau}{3C_f^3M_1^2} \right)\frac{f(g(\bx^0_0)) - F^*}{KJ}                                                      \\
		 & \quad + \left(\frac{A_1C_{h_f}L_f\tau}{2C_f^2M_1} + \frac{2A_3C_{h_f}L_f\tau}{3C_f^3M_1} + \frac{2A_1M_2}{3C_f^2L_gM_1}  + \frac{4A_3M_2}{15C_f^3L_gM_1} + \frac{4A_4M_2}{15C_f^3C_g^2L_gM_1} \right)\epsilon \\
		 & \quad +  \frac{A_2C_{h_f}L_fM_2\tau}{3C_f^3M_1^2} \frac{\left(f(g(\bx^0_0)) - F^*\right)\epsilon}{KJ} +\left( \frac{A_2C_{h_f}L_fM_2\tau}{3C_f^3M_1} + \frac{A_2M_2^2}{18C_f^3L_gM_1} \right)\epsilon^2,
	\end{align*}
	}%
	where the coefficients are defined in the proof.
\end{theorem}

The proof is provided in Appendix \ref{pf:RoS-sample complexity}. By setting $K = J = \cO(\epsilon^{-1/2})$, the  {RoS-VR} algorithm \ref{alg:RoS-SPIDER} achieves its best sample complexity which is $\cO(\epsilon^{-5/2})$ for the $g_\xi$ oracle, $\cO(\epsilon^{-3/2})$ for the $\grad g_\xi$ oracle, and $\cO(\epsilon^{-2})$ for the $\grad f_\varphi$ oracle.



\section{Applications and Numerical Experiments} \label{sec:applications}
 {Below, we evaluate the performance of the proposed algorithms over two different applications and compared them with two other algorithms for stochastic composition optimization. All implementations are performed in Python and are available at \url{https://github.com/samdavanloo/BSC}.}

\subsection{Risk-averse optimization (SoR case) } \label{sec:risk_averse_ex}
 {The mean-variance risk-averse optimization problem \eqref{eq:mean_var} can be written as}
\vspace{-0.2cm}
\begin{equation}\label{eq:risk-averse}
	\min_{\bx \in \cX} 
	-\mathbb{E}\left[r_\xi(\bx)\right] + \lambda \left( \mathbb{E}\left[\left(r_\xi(\bx)\right)^2\right] - \mathbb{E}^2\left[r_\xi(\bx)\right] \right),
\end{equation}
 {which is an instance of two-level stochastic composition problem \eqref{eq:2-level composition} with}
\vspace{-0.2cm}
\begin{equation*}
	g(\bx):\mathbb{R}^n\rightarrow \mathbb{R}^2 \triangleq \mathbb{E}\left[r_\xi(\bx),  \quad \left(r_\xi(\bx)\right)^2 \right], \quad  f(u_1,u_2): \mathbb{R}^2 \rightarrow \mathbb{R} \triangleq -u_1 + \lambda u_2 - \lambda u_1^2.
\end{equation*}
 {Following our discussion in Example~\ref{ex:risk_averse}, we let the random reward to be a quadratic function} $r_\xi(\bx) \triangleq \frac{1}{2} \bx^\intercal A_\xi \bx$ where $A_\xi$ is a symmetric matrix. Furthermore, we let the feasible set to be $\cX \triangleq \{\bx \in \mathbb{R}|\norm{x} \leq R \}$ for some $R>0$.  {Note that while $f$ is quadratic and hence smooth,} $g(\bx)$ does not have Lipschitz continuous Jacobian as $g_{\xi}(\bx)$ involves fourth degree polynomial $(r_\xi(\bx))^2$. However, the following lemma shows that $(r_\xi(\bx))^2$ is relatively smooth which we use in Lemma~\ref{lem:rel_app_1} to show that $g_\xi$ is average relatively smooth.  {Hence, the objective of \eqref{eq:risk-averse} is an instance of the SoR composition.}
\begin{lemma}[Lemma 5.1 in \cite{bolte2018first}]\label{lemma:4th rel-smooth}
	For any $L\geq 3\|A\|^2$, $f(\bx) =  \frac{1}{4}(\bx^\intercal A\bx)^2$ is L-smooth relative to $h(x) = \frac{1}{4}\|\bx\|^4$.
\end{lemma}
 
\begin{lemma} \label{lem:rel_app_1}
	The vector-valued function $g_\xi(\bx) = [ \frac{1}{2}\bx^\intercal A_\xi \bx  \quad \frac{1}{4}(\bx^\intercal A_\xi\bx)^2]$  is average 1-smooth relative to $h_g(\bx) = \frac{k_1}{2}\|\bx\|^2 + \frac{k_2}{4}\|\bx\|^4$, where $k_1 \geq \norm{\mathbb{E}\left[A_\xi\right]}$ and $k_2 \geq 3\mathbb{E}\left[\|A_\xi\|^2\right]$.
\end{lemma}

\begin{remark}\label{remark:SoR_bdd_lip_smooth_}
	Note that when the feasible set $\cX$ is compact, the inner function $g$ is also smooth, but its smoothness constant depends on the diameter of $\cX$ since $\norm{\grad^2 g_2(\bx)} \leq \sup_{\bx\in\cX} 3\mathbb{E}[\|A_\xi\|^2]\norm{\bx}^2 = \cO(R^2)$.
\end{remark}
For the boundedness of the stochastic gradient of $g_\xi$ in expectation, we have 
$\grad g_\xi(\bx) = \begin{bmatrix}
		A_\xi\bx \quad  \bx^\intercal A_\xi \bx A_\xi \bx
	\end{bmatrix}$, so
 {\small
\begin{align*}
	\mathbb{E}[\norm{\grad g_\xi(\bx)}^2]  \leq \mathbb{E}[\norm{\grad g_\xi(\bx)}_F^2]  = \mathbb{E}[\norm{A_\xi\bx}^2 +\norm{\bx^\intercal A_\xi \bx A_\xi \bx }^2]
	                                                                                     \leq \mathbb{E}[\norm{A_\xi}^2]R^2 + \mathbb{E}[\norm{A_\xi}^4]R^6.
\end{align*}
}
Hence, $C_g^2 \geq \mathbb{E}[\norm{A_\xi}^2]R^2 + \mathbb{E}[\norm{A_\xi}^4]R^6$. Furthermore, the gradient and Hessian of $f$ are
\vspace{-0.2cm}
\begin{align*}
	\grad f(u_1,u_2) = \begin{bmatrix}
		-1 - 2\lambda u_1 \\
		\lambda
	\end{bmatrix}, \quad \grad^2 f(u_1,u_2) = \begin{bmatrix}
		-2\lambda & 0 \\
		0         & 0
	\end{bmatrix}.
\end{align*}
\vspace{-0.2cm}
With $A \triangleq \mathbb{E}[A_\xi]$, we have
\vspace{-0.1cm}
\begin{align*}
	\norm{\grad f(u_1,u_2)}^2  & = (1+2\lambda u_1)^2 + \lambda^2 \leq 1+ \lambda^2\norm{A}^2\norm{\bx}^4 + 2\lambda \norm{A}\norm{\bx}^2 + \lambda^2 \\
	                           & \leq 1+\lambda^2R^4\norm{A}^2 + 2\lambda R^2 \norm{A} + \lambda^2,
\end{align*}
and $\norm{\grad ^2 f(u_1,u_2)} = 2\lambda$. So $C_f^2 \geq 1+\lambda^2R^4\norm{A}^2 + 2\lambda R^2 \norm{A} + \lambda^2$. Finally, the average smoothness of $f_\varphi$ and average relative smoothness of $g_\xi$ follow with $L_f = 2\lambda$ and $L_g=1$. From the discussion above and Lemma~\ref{lemma:SoR smooth}, we have $f(g(\bx))$ is 1-smooth relative to $h(\bx) = \frac{C_g^2L_f+C_fL_g\norm{A}}{2}\norm{\bx}^2+ \frac{3 C_fL_g\mathbb{E}[\|A_\xi\|^2] }{4}\norm{\bx}^4$. Note that the above discussion provides loose bounds on $C_f,C_g$. In the experiments, we use the distance generating function for the composition as
\begin{align*}
	h(\bx) \triangleq \frac{ {c_1} }{2}\norm{\bx}^2 + \frac{ {c_2}}{4} \norm{\bx}^4,\numberthis \label{eq:SoR generating function}
\end{align*}
with $ {c_1,c_2}$ determined by grid search. To do so, we generate grids  {$\cC_1, \cC_2$ with 6 logarithmically spaced values in $[10^{-2}, 10^2]$ for $\cC_1$ and in $[10^{-3}, 10]$ for $\cC_2$}. Next, for all combinations $(c_1,c_2)\in \cC_1\times\cC_2$, we run the algorithm for a large enough number of iterations to find the best combination $(c_1^*,c_2^*)$, and, hence, the distance generating function.

In this experiment, we consider solving \eqref{eq:risk-averse} in a finite-sum form, i.e., the expectations are available as finite-sums divided by the number of components. To simulate the data, we first randomly generate a $50 \times 50$ symmetric matrix $A$. Then, we generate 1000 random noise matrices of size $50\times 50$ with independent elements from standard normal distribution and add them to $A$ to construct 1000 samples of $A_\xi$. We solve \eqref{eq:risk-averse} with $\lambda$ equal to  {10} and the radius of the feasible set $R$ equal to 10.

 Following Corollary~\ref{cor:SoR-rate of convergence}, the step sizes are set to $\tau_k \equiv   {0.025}$ and $\beta_k \equiv  {0.5}$. The grid search over the coefficients of \eqref{eq:SoR generating function} results in $(c_1^*,c_2^*) =  {(15.8, 6.3\times 10^{-3})}$.

To validate the sample complexity of our algorithm, we run the experiment with three different batch sizes $\abs{\cB^k_{\grad g}} = \abs{\cB^k_{\grad f}} \in \{1,10,100\}$. Note that the number of calls to the inner and outer function stochastic gradient oracles are set to be equal in Corollary~\ref{cor:SoR-rate of convergence}  {(which we named $|\cB_\grad|$)}. For each setting, we replicate the algorithm 20 times from a fixed initial point. 

\begin{figure}[hbt]
    \centering
    \includegraphics[width=\textwidth]{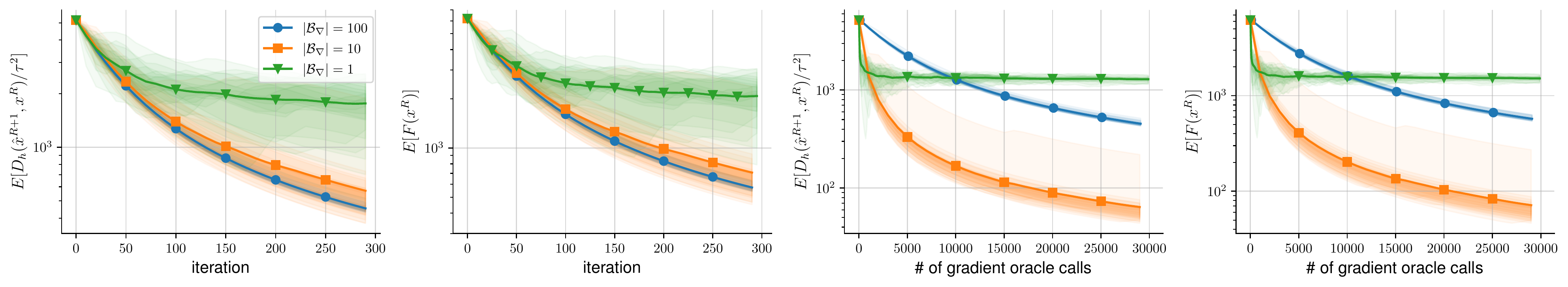}
	\caption{ {Decay of the stationarity measure $\mathbb{E}\left[D_h(\hat{\bx}^{R+1},\bx^R)\right]$ and expected function value $\mathbb{E}[F(\bx^R)]$ versus iteration and stochastic gradient oracle calls by the SoR algorithm~\ref{alg:SoR} for the risk-averse optimization problem.}} 
 \label{fig:SoR}
\end{figure}
Figure \ref{fig:SoR} shows the decay of the stationarity measure  {and function value} versus the iteration and the number of stochastic gradient calls where the solid lines are the average of the 20 replicates. 

Note that based on Corollary~\ref{cor:SoR-rate of convergence} the theoretical sample complexity for $\mathbb{E}[D_h(\hat{\bx}^{R+1},\bx^R)]$ is $\cO(K^{-1}+ \abs{\cB_\grad}^{-1})$. So, for a fixed iteration, the upper bound on stationarity measure is determined by the batch size inverse as verified by the first two plots of Figure~\ref{fig:SoR}.
Comparing the decay of the stationarity measure from sample complexity perspective, we see that $\abs{\cB_\grad} = 10$ performs the best followed by  {$\abs{\cB_\grad} = 100$}. We should also note that with $\abs{\cB_\grad} = 1$ the algorithm has significantly higher variance compared to the bigger two batch sizes. Finally, with $\abs{\cB_\grad} = 10$ the algorithm is able to obtain lower values of the stationarity measure compared to the other two batch sizes.

\begin{figure}[hbt]
    \centering
    \includegraphics[width=\textwidth]{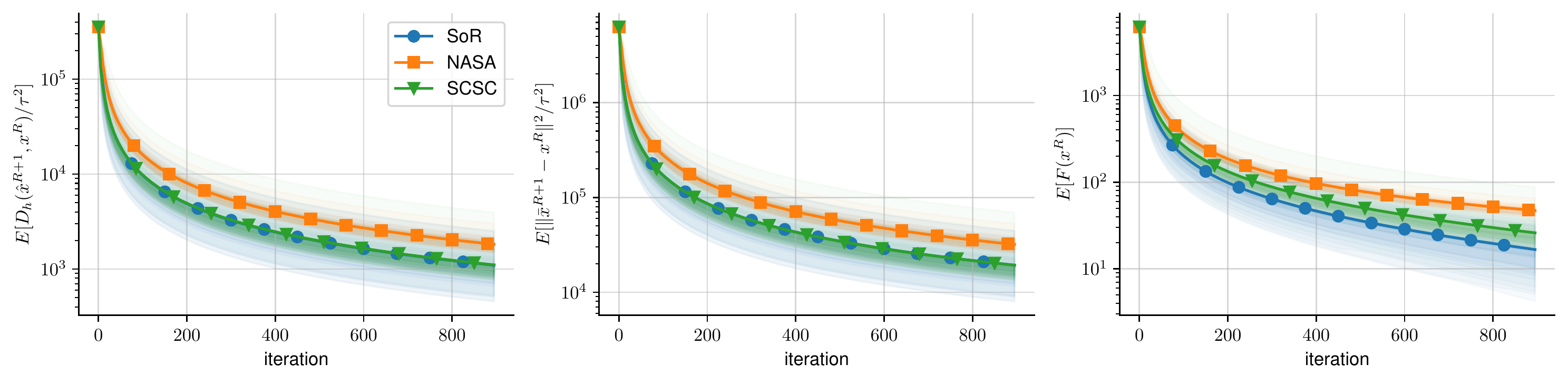}
    \caption{ {Decay of the two stationarity measures and expected function values by the SoR algorithm~\ref{alg:SoR}, NASA~\cite{ghadimi2020single} and SCSC~\cite{chen2020solving} for the risk-averse problem. 
    }}
    \label{fig:SoR_alg_compare}
\end{figure}

We also compare the SoR algorithm with NASA~\cite{ghadimi2020single} and SCSC~\cite{chen2020solving} which are two of the state of the art algorithms for stochastic composition optimization. The batch sizes of the SCSC algorithm are set to 100 and all batch sizes of the NASA algorithm are set to 1 which is supported by their theory. For parameter turning, we also do a grid search and determine the best setting based on the decay in the stationarity measure. 
Figure~\ref{fig:SoR_alg_compare} illustrates the decay of two stationarity measures and function values by different methods. The second stationarity measure is based on $\bar{x}^{k+1}$ defined as
\vspace{-0.2cm}
\begin{align*}
\bar{x}^{k+1} \triangleq \argmin_{\by \in \cX} \fprod{\grad F(\bx^k),\by-\bx^k} + \frac{1}{2\tau_k}\|\by-\bx\|^2,
\end{align*}
which is the Euclidean counterpart of $\hat{x}^{k+1}$ defined in \eqref{eq:opt_obj_1}. SoR and SCSC have similar performance while NASA decays slower across all three measures.

\subsection{Policy Evaluation for MDP (RoS case)}\label{sec:MDP_experiment}
 {Following the discussion in Example~\ref{ex:MDP}, in this section, we solve \eqref{eq:policy evaluation problem} when the reward is of count type (e.g., number of clicks). It is common in this setting to take the distance function to be the KL divergence $D_{KL}(a,b) = \sum_i a_i\log\frac{a_i}{b_i} + b_i - a_i$ with $a_i,b_i >0$ which better captures the information when the random noise follows the Poisson distribution.}

The problem is formulated in \eqref{eq:policy evaluation problem} which is an instance of the stochastic composition problem~\eqref{eq:2-level composition} with
\vspace{-0.2cm}
\begin{align}\label{eq:f_and_g_in_MDP}
	g(\bx) = \mathbb{E}[A_\xi]\bx, \qquad \quad f(\bu) = \sum_{i=1}^{s} \bu_i -\mathbb{E}[ \br_\varphi]_i
	\log\bu_i,
 \vspace{-0.2cm}
\end{align}
where $s\triangleq|\cY|$. Note that in \eqref{eq:f_and_g_in_MDP}, $g_\xi$ is average Lipschitz smooth
and $f_\varphi$ is average smooth relative to $h_f(\bu) = -\sum_{i=1}^{s}\log \bu_i$; hence, this is an instance of the RoS composition. As RoS is a subset of the RoR composition, in this experiment, we solve \eqref{eq:policy evaluation problem} by both the RoR~\ref{alg:RoR-mini batch} and RoS-VR~\ref{alg:RoS-SPIDER} algorithms.

Similar to the SoR experiment, we consider a finite sum setting. To simulate data, we first generate a random matrix $A$ with independent elements $A_{ij} \sim U[0,2],\ i=1,\cdots 50,\ j=1,\cdots,30$, random vector $\br$ with independent elements $\br_i \sim Poisson(0.5)$. Next, we generate 30,000 random noise matrices with independent elements from standard normal distribution add them to $A$ and set the negative elements to zero to create $A_\xi^+$.

The parameters $L_f/\tau$ and $\lambda/\tau$ in \eqref{eq:RoS/R subproblem} and \eqref{eq:RoS subproblem} are  determined by grid search with a grid of 6 logarithmically spaced values in $[10^{-3}, 10^3$]. The RoR algorithm~\ref{alg:RoR-mini batch} is run with two different batch sizes $\abs{\cB_\grad}=100, \abs{\cB_g}=10^4$ and $\abs{\cB_\grad}=20, \abs{\cB_g}=400$. Similarly, the  {RoS-VR} algorithm~\ref{alg:RoS-SPIDER} is run with two different batch sizes $\abs{\cB_\grad}=100, \abs{\cB_g}=10^4,\abs{\cS_\grad}=10, \abs{\cS_g}=100$ and $\abs{\cB_\grad}=20, \abs{\cB_g}=400,\abs{\cS_\grad}=5, \abs{\cS_g}=25$ (the bigger batch sizes for $j = 0$ match the setting of the RoR algorithm and the smaller batch sizes for $j >0$ are \textit{roughly} the square root of the bigger batch sizes).
In the  {RoS-VR} algorithm, the maximum iteration numbers in the two nested loops, i.e., $K$ and $J$, are set equal to each other and equal to the square root of the maximum iteration number in the RoR algorithm. Both algorithms are replicated 20 times from a fixed initial point for each scenario.

\begin{figure}[hbt]
    \centering
    \includegraphics[width = \textwidth]{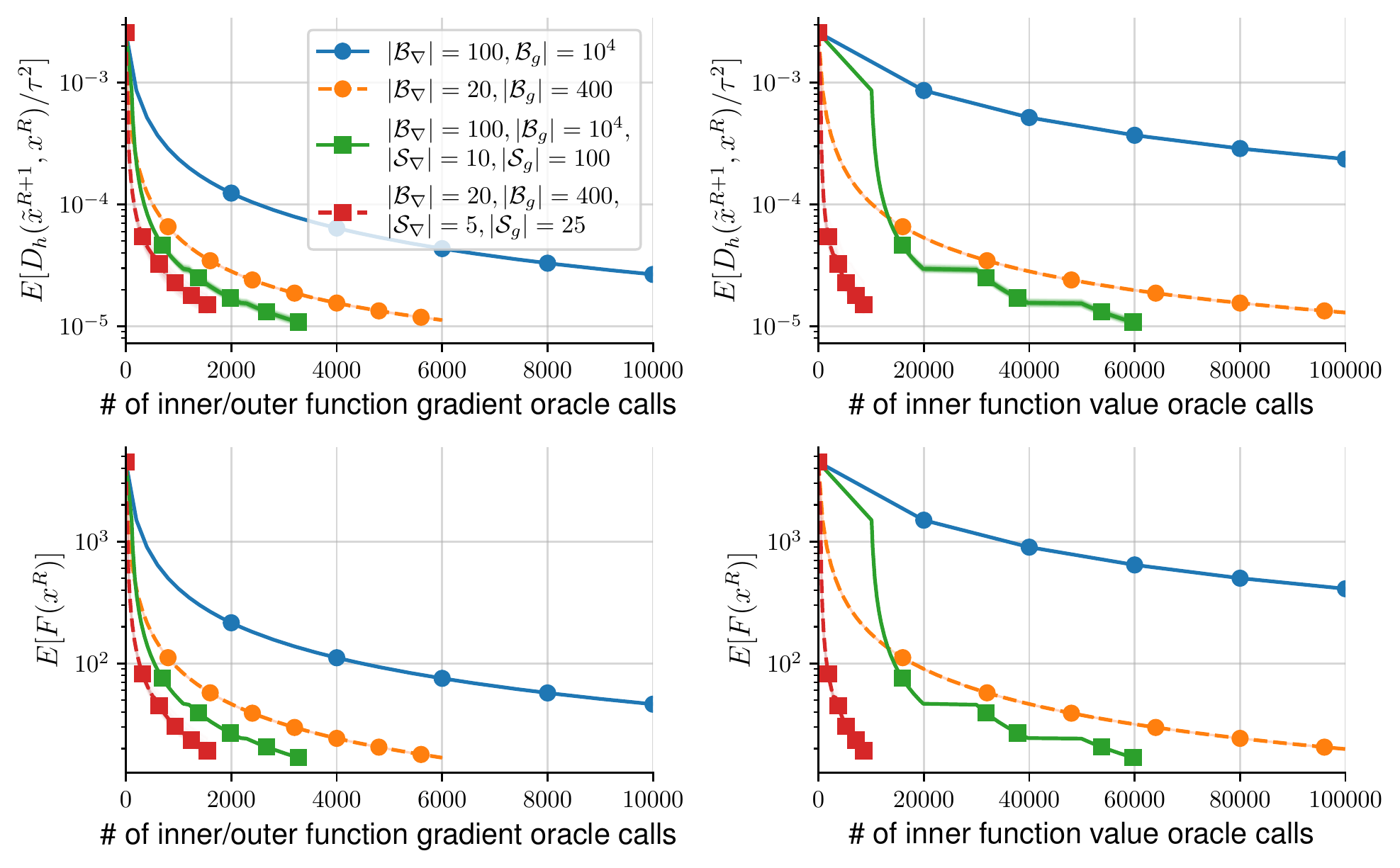}
    \caption{ {Decay of the stationarity measure $\mathbb{E}\left[D_h(\tilde{\bx}^{R+1},\bx^R)\right]$ and expected function value $\mathbb{E}[F(\bx^R)]$ versus the inner/outer gradient and inner function value stochastic oracle calls by the RoR~\ref{alg:RoR-mini batch} (circle marker) and RoS-VR~\ref{alg:RoS-SPIDER} (square marker) algorithms for the policy evaluation problem. 
    Bigger and smaller batch sizes are shown with solid and dashed lines, respectively.}}
	\label{fig:RoS}
\end{figure}

 Figure~\ref{fig:RoS} shows the decay of the stationarity measure and expected function value versus the inner/outer gradient and inner function value stochastic oracle calls by the RoR~\ref{alg:RoR-mini batch} and RoS-VR~\ref{alg:RoS-SPIDER} algorithms over 20 replicates. Results show that the decays are faster by the variance-reduced RoS-VR algorithm compared to the RoR algorithm which supports our theoretical findings.

\begin{figure}
    \centering
    \includegraphics[width = \textwidth]{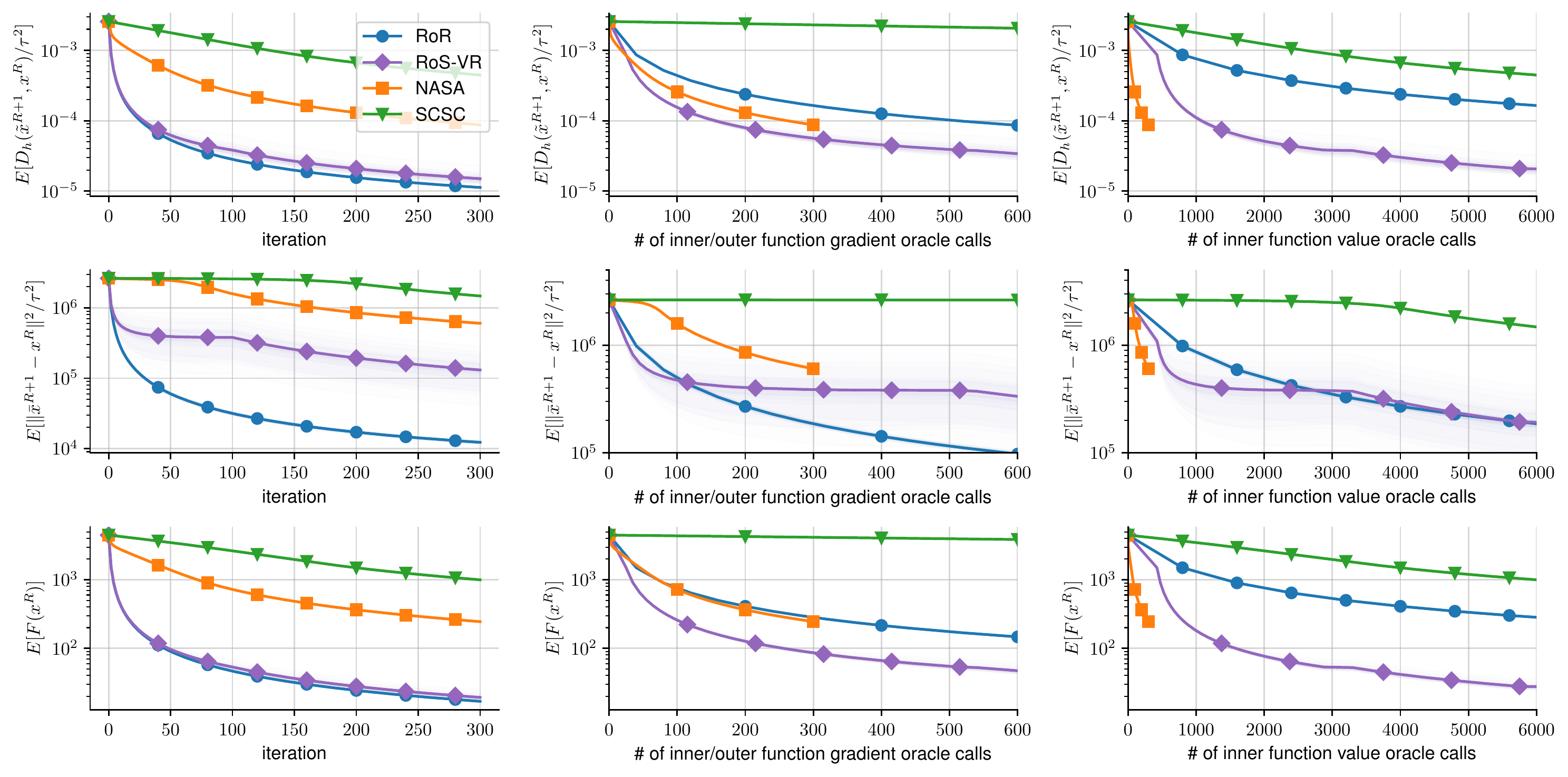}
    \caption{ {Decay of the stationarity measures and expected function values by RoR~\ref{alg:RoR-mini batch}, RoS-VR~\ref{alg:RoS-SPIDER}, NASA~\cite{ghadimi2020single} and SCSC~\cite{chen2020solving} for the policy evaluation problem.}}
    \label{fig:RoS_alg_compare}
\end{figure}

We also compare the two algorithms with NASA~\cite{ghadimi2020single} and SCSC~\cite{chen2020solving}. The results are shown in Figure~\ref{fig:RoS_alg_compare}. The proposed algorithms perform generally better than the other two methods with respect to the iteration and the inner/outer gradient stochastic oracle calls. However, NASA decays faster with respect to the inner function stochastic oracle calls which is expected given its single-batch oracle call per iteration. 
With respect to iteration, RoR has a good performance relative to RoS-VR. But, with respect to the number of stochastic oracle calls, RoS-VR results in faster decay across different measures. For most of the settings, SCSC could not obtain low values of the stationarity measure.

\section{Conclusions}
In this paper, we study the two-level stochastic composition problem in the absence of Lipschitz continuity of the gradient of the inner, outer, or both functions. Given the notion of  relative smoothness for the single-level deterministic problems, we consider three compositions: smooth of relative-smooth (SoR), relative-smooth of smooth (RoS), and relative-smooth of relative-smooth (RoR) compositions.  We then propose one algorithm to solve the SoR and one algorithm to solve the RoR and RoS compositions. We further improve the second algorithm by variance reduction for the RoS composition. We then investigate the iteration and sample complexities of the three proposed algorithms. Finally, we evaluate the performances of these algorithms over two numerical experiments and compare them with other methods.


\vspace{0.2cm}
\bibliographystyle{spmpsci}
\bibliography{refs.bib}


\appendix  
\begin{appendices}



\vspace{-0.5cm}
\section{SoR Proofs}

\vspace{-0.2cm}
\subsection{Proof of Lemma \ref{lemma:SoR smooth}}\label{pf:SoR smooth}
By $L_f$-smoothness of $f$, $\forall g(\bx),g(\by) \in \dom f$, we have
    \begin{align*}
    & \quad \frac{L_f}{2}\norm{g(\bx)-g(\by)}^2 \\
	& \geq \abs{f(g(\bx)) - f(g(\by)) - \fprod{\grad f(g(\by)),g(\bx)-g(\by)} } \\
    & \geq |f(g(\bx)) - f(g(\by)) - \fprod{\grad g(\by)\grad f(g(\by)),\bx-\by}| -\norm{\grad f(g(\by))} \norm{g(\bx)-g(\by)-\grad g(\by)^\intercal(\bx-\by)} \\
	& \geq |f(g(\bx)) - f(g(\by)) - \fprod{\grad g(\by)\grad f(g(\by)),\bx-\by}| - C_f L_g D_{h_g}(\bx,\by),
\end{align*}
where the last inequality holds because $f$ is $C_f$-Lipschitz continuous and $g$ is $L_g$-smooth relative to $h_g$.
Rearranging the above result and using the property that $g$ is $C_g$-Lipschitz continuous, with $h(\bx) \triangleq  \frac{C_g^2L_f }{2}\|\bx\|^2 + C_fL_g h_g(\bx)$, we have
\begin{equation*}
	|f(g(\bx)) - f(g(\by)) - \fprod{\grad g(\by)\grad f(g(\by)),\bx-\by}| \leq D_h(\bx,\by).
\end{equation*}


\vspace{-0.5cm}
\subsection{Proof of Lemma \ref{lemma:SoR-error of u}}\label{pf:SoR-error of u}
Define $\bar{g}_\xi(\cdot)\triangleq  \frac{1}{\abs{\cB_g^{ {k+1}}}}\sum_{\xi \in \cB_g^{ {k+1}}}g_{\xi}(\cdot)$, we have
$
	\mE[\|{\bar{g}_\xi(\bx^{k+1}) - g(\bx^{k+1})}\|^2|\cF^{k}] \leq \sigma_g^2/\abs{\cB^{ {k+1}}_g}
$.
Furthermore, by \eqref{eq:u update}, we have
\begin{align*}
	 & \quad \mE[\|g(\bx^{k+1}) - \bu^{k+1}\|^2 | \cF_{k}]  \\
	 & =\mE[\|(1-\beta_k) (g(\bx^k) - \bu^k) + (1 - \beta_k) (g(\bx^{k+1}) - g(\bx^k)) + \beta_k (g(\bx^{k+1}) - \bar{g}_\xi(\bx^{k+1}) )  \\
	 & \quad +(1-\beta_k)( \bar{g}_\xi(\bx^k) - \bar{g}_\xi(\bx^{k+1}) )\|^2 | \cF_{k}]  \\
	 & \leq (1-\beta_k)^2\|g(\bx^k) - \bu^k\|^2  +  \mE[ 2\|(1-\beta_k)(g(\bx^{k+1})-g(\bx^k)) + \beta_k (g(\bx^{k+1}) - \bar{g}_\xi(\bx^{k+1}))\|^2  \\
	 & \quad + 2\|(1-\beta_k)( \bar{g}_\xi(\bx^k) - \bar{g}_\xi(\bx^{k+1}) )\|^2 | \cF_{k}]  \\
	 & = (1-\beta_k)^2\|g(\bx^k) - \bu^k\|^2 + 2(1-\beta_k)^2\|g(\bx^{k+1}) - g(\bx^k)\|^2 + 2\beta_k^2 \mE[ \|g(\bx^{k+1}) - \bar{g}_\xi(\bx^{k+1})\|^2|\cF_{k}] \\
	 & \quad + 4(1-\beta_k)\beta_k \mE[\langle g(\bx^{k+1}) - g(\bx^k)),g(\bx^{k+1}) - \bar{g}_\xi(\bx^{k+1})\rangle|\cF_{k}]  \\
	 & \quad + 2(1-\beta_k)^2 \mE[\|\bar{g}_\xi(\bx^k) - \bar{g}_\xi(\bx^{k+1}) \|^2 | \cF_{k}]  \\
	 & \stackrel{(a)}{\leq} (1-\beta_k)^2 \|g(\bx^k) - \bu^k\|^2 + 4(1-\beta_k)^2 C_g^2\|\bx^{k+1} - \bx^k\|^2 + {2\beta_k^2\sigma_g^2}/{\abs{\cB_g^{ {k+1}}}},
\end{align*}
where $(a)$ is based on the Lipschitz continuity of $g$.


\vspace{-0.5cm}
\subsection{Proof of Lemma \ref{lemma:SoR-error of w} }\label{pf:SoR-error of w}
By Algorithm \ref{alg:SoR},  we have
\begin{align*}
	 & \quad \|\grad F(\bx^k) - \bw^k\|^2  \\
	 & = \|\grad g(\bx^k)\grad f(g(\bx^k)) - \grad g(\bx^k)\bs^k + \grad g(\bx^k)\bs^k - \bv^k\bs^k\|^2  \\
	 & \leq 2\|\grad g(\bx^k)\|^2\|\grad f(g(\bx^k)) - \bs^k\|^2 + 2\|\grad g(\bx^k) - \bv^k\|^2\|\bs^k\|^2  \\
	 & \leq 2C_g^2 \|\grad f(g(\bx^k))- \grad f(\bu^k) + \grad f(\bu^k) - \bs^k\|^2 + 2\|\grad g(\bx^k) - \bv^k\|^2\|\bs^k\|^2  \\
	 & \leq 2C_g^2L_f^2\|g(\bx^k)- \bu^k\|^2+ 2C_g^2\|\grad f(\bu^k) - \bs^k\|^2  \\
	 & \quad + 4C_g^2 \fprod{\grad f(g(\bx^k)) - \grad f(\bu^k),\grad f(\bu^k) - \bs^k}+ 2\|\grad g(\bx^k) - \bv^k\|^2\|\bs^k\|^2,
\end{align*}
where the first inequality follows from the identity $(a+b)^2\leq 2a^2+2b^2$, the second inequality uses the Lipschitz continuity of $g$, and the last inequality uses the smoothness of $f$. Taking  expectation on both sides and using the tower property in the three identities
\begin{align*}
	\mE[\|\grad f(\bu^k) - \bs^k\|^2 ] = \mE_{\bu^k} \left[\mE[\|\grad f(\bu^k) - \bs^k\|^2 |\bu^k ] \right] \leq {C_f^2}/{|\cB^{ {k}}_{\grad f}|},
\end{align*}
\begin{align*}
	 & \quad\mE[\langle\grad f(g(\bx^k)) - \grad f(\bu^k),\grad f(\bu^k) - \bs^k\rangle]
	\\
	 & = \mE_{\bu^k}\left[\mE[ \langle \grad f(g(\bx^k)) - \grad f(\bu^k),\grad f(\bu^k) - \bs^k\rangle|\bu^k]\right] = 0,
\end{align*}
and
\begin{align*}
	 \mathbb{E}[\|\grad g(\bx^k) - \bv^k\|^2\|\bs^k\|^2 ]                                                                   
	 & = \mathbb{E}_{\bu^k}\left[\mathbb{E}[\|\grad g(\bx^k) - \bv^k\|^2\|\bs^k\|^2 |\bu^k] \right]  \\
	 & =  \mathbb{E}_{\bu^k}\left[\mathbb{E}[\|\grad g(\bx^k) - \bv^k\|^2 | \bu^k] \mathbb{E}[\|\bs^k\|^2 |\bu^k] \right]  \leq {C_g^2C_f^2}/{|\cB_{\grad g}^{ {k}}|},
\end{align*}
completes the proof.


\vspace{-0.5cm}
\subsection{Proof of Lemma \ref{lemma:SoR-Vk decrease}}\label{pf:SoR-Vk decrease}
By Definition~\ref{def:rel_smooth} and Lemma \ref{lemma:SoR smooth}, $\forall \bx,\by,\bz \in \cX$, we have
\begin{align*}
	F(\bx)- F(\by)  & \leq \fprod{\grad F(\bx) , \bx - \by} + D_h(\by,\bx), \\
	F(\bz) - F(\bx) & \leq \fprod{\grad F(\bx), \bz - \bx} + D_h(\bz,\bx).
\end{align*}
Summing the above two inequalities, we get
\begin{equation}\label{eq:smooth ineq}
	F(\bz) - F(\by) \leq \fprod{\grad F(\bx),\bz - \by}+ D_h(\by,\bx) + D_h(\bz,\bx).
\end{equation}
By {Lemma~\ref{lemma:three-point ineq}}, $\forall \by\in\cX$, we have
\begin{align}
	\tau_k\langle\bw^k,\bx^{k+1}-\by\rangle \leq 	D_h(\by,\bx^k)-D_h(\bx^{k+1},\bx^k) - D_h(\by,\bx^{k+1}), \label{eq:three-point of w} \\
	\frac{\tau_k}{2}\langle\grad F(\bx^k),\hat{\bx}^{k+1}-\by\rangle \leq 	D_h(\by,\bx^k)-D_h(\hat{\bx}^{k+1},\bx^k) - D_h(\by,\hat{\bx}^{k+1}). \label{eq:three-point of grad}
\end{align}
Letting $\bz = \bx^{k+1},\bx = \bx^k, \by = \hat{\bx}^{k+1}$ in  \eqref{eq:smooth ineq} and \eqref{eq:three-point of w}, and summing them up, we have
\begin{align*}
	F(\bx^{k+1}) - F(\hat{\bx}^{k+1})
	 & \leq \langle\grad F(\bx^k)-\bw^k,\bx^{k+1} - \hat{\bx}^{k+1}\rangle + \left(1+\frac{1}{\tau_k}\right)D_h(\hat{\bx}^{k+1},\bx^k) \\
	 & \quad   + \left(1 - \frac{1}{\tau_k}\right)D_h(\bx^{k+1},\bx^k)  - \frac{1}{\tau_k}D_h(\hat{\bx}^{k+1},\bx^{k+1}).
\end{align*}
Letting $\bz = \hat{\bx}^{k+1}, \bx = \bx^k, \by = \bx^k$ in \eqref{eq:smooth ineq} and \eqref{eq:three-point of grad}, summing them up,
\begin{equation*}
	F(\hat{\bx}^{k+1}) - F(\bx^k) \leq \left(1 - \frac{2}{\tau_k}\right)D_h(\hat{\bx}^{k+1},\bx^k) - \frac{2}{\tau_k}D_h(\bx^k,\hat{\bx}^{k+1}).
\end{equation*}
Adding the last two inequalities, we get
\begin{align*}
	 & \quad  F(\bx^{k+1}) - F(\bx^k)  \\
	 & \leq \langle\grad F(\bx^k)-\bw^k,\bx^{k+1} - \hat{\bx}^{k+1}\rangle + \left(2 - \frac{1}{\tau_k} \right)D_h(\hat{\bx}^{k+1},\bx^k)   + \left(1 - \frac{1}{\tau_k}\right)D_h(\bx^{k+1},\bx^k)  \\
	 & \quad - \frac{1}{\tau_k}D_h(\hat{\bx}^{k+1},\bx^{k+1}) - \frac{2}{\tau_k}D_h(\bx^k,\hat{\bx}^{k+1})  \\
	 & \leq \frac{\tau_k}{2\mu}\|\grad F(\bx^k) - \bw^k\|^2 + \frac{\mu}{2\tau_k}\|\bx^{k+1} - \hat{\bx}^{k+1}\|^2 + \left(2 - \frac{1}{\tau_k} \right)D_h(\hat{\bx}^{k+1},\bx^k)  \\
	 & \quad  + \left(1 - \frac{1}{\tau_k}\right)D_h(\bx^{k+1},\bx^k)     - \frac{1}{\tau_k}D_h(\hat{\bx}^{k+1},\bx^{k+1})  \\
	 & \leq \frac{\tau_k}{2\mu}\|\grad F(\bx^k) - \bw^k\|^2 + \frac{1}{\tau_k}D_h(\hat{\bx}^{k+1},\bx^{k+1})+ \left(2 - \frac{1}{\tau_k} \right)D_h(\hat{\bx}^{k+1},\bx^k)  \\
	 & \quad + \left(1 - \frac{1}{\tau_k}\right)D_h(\bx^{k+1},\bx^k) - \frac{1}{\tau_k}D_h(\hat{\bx}^{k+1},\bx^{k+1})  \\
	 & \leq \frac{\tau_k}{2\mu}\|\grad F(\bx^k) - \bw^k\|^2  + \left(2 - \frac{1}{\tau_k} \right)D_h(\hat{\bx}^{k+1},\bx^k)   +  \left(1 - \frac{1}{\tau_k}\right)D_h(\bx^{k+1},\bx^k),                                                                               \numberthis \label{eq:F_ub}
\end{align*}
where the second inequality is due to Young's inequality and the third inequality holds because $h$ is $\mu$-strongly convex.

Adding $-F(\bx^*)+\|g(\bx^{k+1})-\bu^{k+1}\|^2$ on both sides of the inequality \eqref{eq:F_ub}, taking expectation conditioned on $\cF_k$, we get
\begin{align*}
	 & \quad \mE[F(\bx^{k+1}) - F(\bx^*)+ \|g(\bx^{k+1})-\bu^{k+1}\|^2|\cF_k]  \\
	 & \leq  F(\bx^k) - F(\bx^*) +\|g(\bx^{k})-\bu^{k}\|^2 + \frac{\tau_k}{2\mu}\|\grad F(\bx^k) -\bw^k\|^2- \left(\frac{1}{\tau_k}-2\right)D_h(\hat{\bx}^{k+1},\bx^k) \\
	 & \quad - \left(\frac{1}{\tau_k} - 1\right)D_h(\bx^{k+1},\bx^k) +\mE[\|g(\bx^{k+1})-\bu^{k+1}\|^2|\cF_k]-\|g(\bx^{k})-\bu^{k}\|^2.
\end{align*}
Using the definition of $V^k$, Lemmas \ref{lemma:SoR-error of u} and the fact that $\mu = C_g^2L_f$, we have
\begin{align*}
	 & \quad \mE[V^{k+1}|\cF_k] \\
	 & \leq V^{k}  + \frac{\tau_k}{2C_g^2L_f}\|\grad F(\bx^k) -\bw^k\|^2 - \left(\frac{1}{\tau_k} - 2 \right)D_h(\hat{\bx}^{k+1},\bx^k)    - \left(\frac{1}{\tau_k} - 1 \right)D_h(\bx^{k+1},\bx^k)  \\
	 & \quad   + ((1-\beta_k)^2-1)\|g(\bx^k)-\bu^k\|^2   + 4(1-\beta_k)^2C_g^2 \|\bx^{k+1}-\bx^k\|^2 + {2\beta_k^2\sigma_g^2}/{\abs{\cB_g^{ {k+1}}}}  \\
	 & \leq V^k + \frac{\tau_k}{2C_g^2L_f}\|\grad F(\bx^k) -\bw^k\|^2 - \left(\frac{1}{\tau_k}-2\right)D_h(\hat{\bx}^{k+1},\bx^k)  - \left(\frac{1}{\tau_k}-1 - \frac{8}{L_f} \right)D_h(\bx^{k+1},\bx^k) \\
	 & \quad  +((1-\beta_k)^2-1)\|g(\bx^k)-\bu^k\|^2 + {2\beta_k^2\sigma_g^2}/{\abs{\cB_g^{ {k+1}}}},
\end{align*}
where the second inequality holds because  $h$ is $\mu$-strongly convex.


\vspace{-0.5cm}
\subsection{Proof of Theorem \ref{thm:SoR-rate of convergence}}\label{pf:SoR-rate of convergence}
Taking expectation of \eqref{eq:bound} with respect to the random sequences generated by the algorithm and using the tower property, we get
\begin{align*}
	\mE[V^{k+1}] & \leq \mE[V^{k}] + \frac{\tau_k}{2C_g^2L_f}\mE[\|\grad F(\bx^k) -\bw^k\|^2] - \left(\frac{1}{\tau_k} - 2 \right)\mE[D_h(\hat{\bx}^{k+1},\bx^k)]  \\
	                    & \quad - \left(\frac{1}{\tau_k}-1 - \frac{8}{L_f} \right)\mE[D_h(\bx^{k+1},\bx^k)] +((1-\beta_k)^2-1)\mE[\|g(\bx^k)-\bu^k\|^2] + {2\beta_k^2\sigma_g^2}/{\abs{\cB^{ {k+1}}_g}}.
\end{align*}
Applying Lemma \ref{lemma:SoR-error of w} to the above inequality, under the assumption$(1-\beta_k)^2+\tau_kL_f-1 \leq 0$ with $0<\tau_k<\min\{1/2,L_f/(L_f+8)\},\beta_k \in (0,1)$, we have
\begin{align*}
	\mE[V^{k+1}] & \leq \mE[V^{k}] - \frac{1-2\tau_k}{\tau_k}\mE[D_h(\hat{\bx}^{k+1},\bx^k)] + \frac{\tau_k C_f^2}{L_f}\left(1/|\cB^{ {k}}_{\grad f}| + 1/|\cB^{ {k}}_{\grad g}| \right)  \\
	                    & \quad + \left((1-\beta_k)^2+\tau_k L_f-1\right)\mE[ \|g(\bx^k)-\bu^k\|^2]  + {2\beta_k^2\sigma_g^2}/{\abs{\cB_g^{ {k+1}}}}  \\
	                    & \leq  \mE[V^{k}] - \frac{1-2\tau_k}{\tau_k}\mE[D_h(\hat{\bx}^{k+1},\bx^k)] + \frac{\tau_k C_f^2}{L_f}\left(1/|\cB^{ {k}}_{\grad f}| + 1/|\cB^{ {k}}_{\grad g}| \right)   + {2\beta_k^2\sigma_g^2}/{\abs{\cB_g^{ {k+1}}}} .
\end{align*}
Rearranging the terms in the above inequality, telescoping from $k=0,\dots,K-1$ and using the fact that $V^{k+1}\geq 0$, we have
\begin{align*}
	\sum_{k=0}^{K-1}(\tau_k - 2\tau_k^2)\mE\left[{D_h(\hat{\bx}^{k+1},\bx^k)}/{\tau_k^2}\right]
	\leq V^0 + \sum_{k=0}^{K-1} \frac{\tau_k C_f^2}{L_f|{\cB^{ {k}}_{\grad f}}|} + \frac{\tau_k C_f^2}{L_f|{\cB^{ {k}}_{\grad g}}|} + \frac{2\beta_k^2\sigma_g^2}{|{\cB_g^{ {k+1}}}|}  .
\end{align*}
Dividing both sides by $\sum_{j=0}^{K-1}(\tau_j - 2\tau_j^2)$ and using the fact that
\begin{align*}
	\mE\left[ {D_h(\hat{\bx}^{R+1},\bx^R)}/{\tau_R^2} \right] =\frac{\sum_{k=0}^{K-1}(\tau_k- 2\tau_k^2)\mE\left[{D_h(\hat{\bx}^{k+1},\bx)}/{\tau_k^2} \right]}{\sum_{j=0}^{K-1}(\tau_j-2\tau_j^2)} ,
\end{align*}
we get
\begin{align*}
	\mE\left[ {D_h(\hat{\bx}^{R+1},\bx^R)}/{\tau_R^2} \right]
	\leq \frac{V^0}{\sum_{j=0}^{K-1}(\tau_j - 2\tau_j^2)}+\frac{\sum_{k=0}^{K-1} \frac{\tau_k C_f^2}{L_f|\cB^{ {k}}_{\grad f}|} + \frac{\tau_k C_f^2}{L_f|\cB^{ {k}}_{\grad g}|} + \frac{2\beta_k^2\sigma_g^2}{|{\cB_g^{ {k+1}}}|}}{\sum_{j=0}^{K-1}(\tau_j - 2\tau_j^2)} .
\end{align*}

\vspace{-0.5cm}
\subsection{Proof of Corollary \ref{cor:SoR-rate of other estimators}}\label{pf:SoR-rate of other estimators}
By Lemma \ref{lemma:SoR-error of u}, the $C_g^2L_f$-strongly convexity of $h$ and the tower property, denoting $\Delta^k \triangleq \norm{g(\bx^k) - \bu^k}^2$, we have
\begin{align*}
	\mE[\Delta^{k+1}] 
	                          \leq (1-L_f^2\tau^2)\mE[\Delta^k] + \frac{8}{L_f}\mE[D_h(\bx^{k+1},\bx^k)] +\frac{\epsilon}{2}(\tau-2\tau^2).
\end{align*}
Rearranging the above terms,
\begin{align*}
	\mE[\Delta^k]            & \leq \mE\left[\frac{\Delta^k-\Delta^{k+1}}{L_f^2\tau^2} \right] + \frac{8}{L_f^3}\mE\left[{D_h(\bx^{k+1},\bx^k)}/{\tau^2}\right] + \frac{(1-2\tau)\epsilon}{2L_f^2\tau} .
\end{align*}
Telescoping the above inequality from $k=0$ to $K-1$, and dividing by $K$ on both sides, we get
\begin{align*}
	\frac{1}{K}\sum_{k=0}^{K-1}\mE[\Delta^k] & \leq \frac{\Delta^0}{KL_f^2\tau^2} + \frac{8}{L_f^3}\frac{1}{K}\sum_{k=0}^{K-1}\mE\left[{D_h(\bx^{k+1},\bx^k)}/{\tau^2}\right] + \frac{(1-2\tau)\epsilon}{2L_f^2\tau}  \\
	                                                & \leq \frac{\Delta^0}{KL_f^2\tau^2} + \frac{8V^0}{KL_f^3(\tau - \frac{L_f + 8}{L_f}\tau^2)} + \left(\frac{8(1-2\tau)}{L_f^3(1 - \frac{L_f + 8}{L_f}\tau)} + \frac{1-2\tau}{2L_f^2\tau} \right)\epsilon.
\end{align*}
{Inserting the above result into \eqref{eq:gradient variance} from Lemma \ref{lemma:SoR-error of w}, the error of the gradient of the composition function can be bounded.}

\vspace{-0.5cm}
\section{RoR Proofs}
\subsection{Proof of Lemma \ref{lemma:RoR-relative smoothness}}\label{pf:RoR-relative smoothness}
$f$ is $L_f$-smooth relative to $h_f$, so $\forall \bx,\by \in \cX$, we have
\begin{align*}
	 & \quad L_f D_{h_f}(g(\bx),g(\by))  \\
	 & \geq \abs{f(g(\bx))-f(g(\by)) - \fprod{\grad f(g(\by)),g(\bx)-g(\by)}}  \\
	 & \geq \abs{f(g(\bx))-f(g(\by)) -\fprod{\grad g(\by)\grad f(g(\by)),\bx-\by}}                                   -C_f\norm{g(\bx)-g(\by) - \grad g(\by)^\intercal(\bx-\by)} \\
	                           & \geq  \abs{f(g(\bx))-f(g(\by)) -\fprod{\grad g(\by)\grad f(g(\by)),\bx-\by}} - C_fL_gD_{h_g}(\bx,\by),
\end{align*}
where the fist inequality follows from the reverse triangle inequality, the second inequality follows from Cauchy-Schwarz inequality, and the last inequality is due to $C_f$-Lipschitz continuity of $f$ and $L_g$-relative smoothness of $g$. Next, we upper bound $D_{h_f}(g(\bx) ,g(\by))$ as
\begin{align*}
	D_{h_f}(g(\bx),g(\by))
	 & =h_f(g(\bx))-h_f(g(\by))  -\langle \grad h_f(g(\by)),g(\bx)-g(\by) \rangle  \\
	 & \quad  +\langle\grad g(\by)\grad h_f(g(\by)),\bx-\by\rangle - \langle\grad g(\by)\grad h_f(g(\by)),\bx-\by\rangle  \\
	                       & \leq A  + \|\grad h_f(g(\by))\|\| g(\bx) -g(\by) - \grad g(\by)^{\intercal}(\bx-\by)\|                             \\
	                       & \leq A + C_{h_f}L_gD_{h_g}(\bx,\by) , \numberthis
\end{align*}
where $A\triangleq  h_f(g(\bx)) - h_f(g(\by)) - \langle\grad g(\by)\grad h_f(g(\by)),\bx-\by\rangle$. Combining the above two inequalities, we have
\begin{align*}
	 & \quad \abs{f(g(\bx)) - f(g(\by)) - \fprod{\grad g(\by)\grad f(g(\by)),\bx-\by}}  \\
	& \leq L_f A +  L_fC_{h_f}L_gD_{h_g}(\bx,\by) + C_fL_gD_{h_g}(\bx,\by) = D_h(\bx,\by),
\end{align*}
where $h(\bx) \triangleq  (C_fL_g+C_{h_f}L_g L_f) h_g(\bx) + L_f h_f(g(\bx))$. Below, we show that the function $h$ is indeed convex.

Using Proposition~\ref{prop:composition_rwc}, $h_f(g(\bx))$ is $C_{h_f}L_g-$weakly convex relative to $h_g$; hence, $C_{h_f}L_g L_f h_g(\bx)+L_fh_f(g(\bx))$ is convex, and hence $h$, is convex. Furthermore, if $h_g$ is $1$-strongly convex, $h$ is $C_fL_g$-strongly convex.

\vspace{-0.5cm}
\subsection{Proof of Lemma \ref{lemma:RoR-linearization bound}}\label{pf:RoR-linearization bound}
Since $h_f$ is $C_{h_f}$-Lipschitz continuous and $g$ is $L_g$-smooth relative to $h_g$, we have
\begin{align*}
	 & \quad |h_f(g(\by)) - h_f(g(\bx)+ \fprod{\grad g(\bx),\by-\bx})|  \\
	 & \leq C_{h_f}\norm{g(\by) - g(\bx) - \fprod{\grad g(\bx),\by-\bx}} 	\leq C_{h_f} L_g D_{h_g}(\by,\bx).
\end{align*}
Rearranging the above terms, the first result follows.
By Lemma \ref{lemma:RoR-relative smoothness}, we have
\begin{align*}
	 & \quad F(\by) \\
	 & \leq F(\bx)+\fprod{\grad F(\bx),\by-\bx} + D_h(\by,\bx) \\
	 & \leq F(\bx)+\fprod{\grad F(\bx) ,\by-\bx} + (C_fL_g + C_{h_f}L_gL_f) D_{h_g}(\by,\bx)   +C_{h_f}L_gL_f D_{h_g}(\by,\bx)\\
	 & \quad  +L_fh_f(g(\bx) + \fprod{\grad g(\bx),\by-\bx}) - L_f h_f(g(\bx)) - L_f \fprod{\grad g(\bx)\grad h_f(g(\bx)),\by-\bx}  \\
	 & = F(\bx)+\fprod{\grad F(\bx) ,\by-\bx} + (C_fL_g + 2C_{h_f}L_gL_f) D_{h_g}(\by,\bx)  \\
	       & \quad + L_f D_{h_f}(g(\bx) + \fprod{\grad g(\bx),\by-\bx},g(\bx)),
\end{align*}
where the second inequality uses the first result.


\vspace{-0.5cm}
\subsection{Proof of Lemma \ref{lemma:RoR-decrease lemma}}\label{pf:RoR-decrease lemma}
Since $f$ is $C_f$-Lipschitz continuous, we have
\begin{align}\label{eq:ub 3}
	|f(\bu^k) - f(g(\bx^k))| \leq C_f\|\bu^k - g(\bx^k)\|,
\end{align}
and
\begin{align*}
	 & \quad \abs{f(g(\bx^{k+1})) - f\big(\bu^k+(\bv^k)^\intercal(\bx^{k+1}-\bx^k)\big)}  \\
	 & \leq C_f\|g(\bx^{k+1}) - \bu^k - (\bv^k)^\intercal(\bx^{k+1}-\bx^k)\|  \\
	 & \leq C_f\|g(\bx^{k+1}) - g(\bx^k)- \grad g(\bx^k)^\intercal (\bx^{k+1}-\bx^k)\|   + C_f \|g(\bx^k) - \bu^k\|  \\
	 & \quad  + C_f \|(\grad g(\bx^k) - \bv^k)^\intercal(\bx^{k+1}-\bx^k)\|  \\
	 & \leq C_fL_gD_{h_g}(\bx^{k+1},\bx^k) +  C_f \|g(\bx^k) - \bu^k\| + C_f \|\grad g(\bx^k) - \bv^k\| \|\bx^{k+1}-\bx^k\|,\numberthis \label{eq:ub 1}
\end{align*}
where the last inequality uses the fact $g$ is $L_g$-smooth relative to $h_g$. $f$ is also $L_f$-smooth relative to $h_f$, so we have
\begin{align*}
	 & \quad |f\big(\bu^k + (\bv^k)^\intercal(\bx^{k+1}-\bx^k)\big) - f(\bu^k) - \langle\bv^k\bs^k,\bx^{k+1}-\bx^k\rangle|  \\
	 & \quad + \langle\grad f(\bu^k) - \bs^k, (\bv^k)^\intercal(\bx^{k+1} - \bx^k)\rangle|  \\
	 & \leq L_f D_{h_f}(\bu^k + (\bv^k)^\intercal(\bx^{k+1}-\bx^k),\bu^k) + \|\bv^k\|\|\grad f(\bu^k) - \bs^k\| \|\bx^{k+1}-\bx^k\|. \numberthis\label{eq:ub 2}
\end{align*}
Combining these three inequalities, we have
\begin{align*}
	 & \quad f(g(\bx^{k+1}))  \\
	 & \leq f(\bu^k+(\bv^k)^\intercal(\bx^{k+1}-\bx^k)) + C_fL_gD_{h_g}(\bx^{k+1},\bx^k) + C_f \|g(\bx^k) - \bu^k\|  \\
	 & \quad + C_f \|\grad g(\bx^k) - \bv^k\| \|\bx^{k+1}-\bx^k\|  \\
	 & \leq f(g(\bx^k)) + \langle\bv^k\bs^k,\bx^{k+1}-\bx^k\rangle  + \frac{L_f}{\tau_k} D_{h_f}(\bu^k + (\bv^k)^\intercal(\bx^{k+1}-\bx^k),\bu^k)       +\frac{\lambda}{\tau_k}D_{h_g}(\bx^{k+1},\bx^k) \\
	 & \quad  + 2C_f \|g(\bx^k) - \bu^k\|+ C_f \|\grad g(\bx^k) - \bv^k\| \|\bx^{k+1}-\bx^k\|                                       +\|\bv^k\|\|\grad f(\bu^k) - \bs^k\| \|\bx^{k+1}-\bx^k\|  \\
	                & \quad + L_f(1-\frac{1}{\tau_k})D_{h_f}(\bu^k + (\bv^k)^\intercal(\bx^{k+1}-\bx^k),\bu^k) + (C_fL_g - \frac{\lambda}{\tau_k} ) D_{h_g}(\bx^{k+1},\bx^k),  \numberthis \label{eq:ub 4}
\end{align*}
where $\lambda = C_fL_g + 2C_{h_f}L_gL_f$ as defined in Lemma~\ref{lemma:RoR-linearization bound} . By the definition of $\bx^{k+1}$, we know that
\begin{align*}
	 & \quad \langle\bv^k \bs^{k},\bx^{k+1}-\bx^k\rangle  + \frac{L_f}{\tau_k} D_{h_f}(\bu^k + (\bv^k)^\intercal(\bx^{k+1}-\bx^k),\bu^k)+\frac{\lambda}{\tau_k}D_{h_g}(\bx^{k+1},\bx^k) \\
	 & \leq \langle\bv^k\bs^{k},\bx^k-\bx^k\rangle  + \frac{L_f}{\tau_k} D_{h_f}(\bu^k + (\bv^k)^\intercal(\bx^k-\bx^k),\bu^k)+\frac{\lambda}{\tau_k}D_{h_g}(\bx^k,\bx^k) = 0.
\end{align*}
Furthermore, given that $\tau_k \leq 1$, we have $L_f(1-\frac{1}{\tau_k})D_{h_f}(\bu^k + (\bv^k)^\intercal(\bx^{k+1}-\bx^k),\bu^k) \leq 0$. Using these two inequalities in the right hand side of  \eqref{eq:ub 4}, we have
\begin{align*}
	f(g(\bx^{k+1})) & \leq f(g(\bx^k)) + 2C_f \|g(\bx^k) - \bu^k\|+ C_f \|\grad g(\bx^k) - \bv^k\|\|\bx^{k+1}-\bx^k\|  \\
	                & \quad  +\|\bv^k\|\|\grad f(\bu^k) - \bs^k\| \|\bx^{k+1}-\bx^k\| + (C_fL_g - \frac{\lambda}{\tau_k} ) D_{h_g}(\bx^{k+1},\bx^k)  \\
	                & \leq f(g(\bx^k)) + 2C_f\|g(\bx^k)-\bu^k\| + \frac{C_f^2}{2}\|\grad g(\bx^k) - \bv^k\|^2 + \frac{1}{2}\|\bx^{k+1}-\bx^k\|^2  \\
	                & \quad + \frac{1}{2}\|\bv^k\|^2 \|\grad f(\bu^k)-\bs^k\|^2 + \frac{1}{2}\|\bx^{k+1}-\bx^k\|^2 + (C_fL_g - \frac{\lambda}{\tau_k} ) D_{h_g}(\bx^{k+1},\bx^k) \\
	                & \leq f(g(\bx^k)) + 2C_f\|g(\bx^k)-\bu^k\| + \frac{C_f^2}{2}\|\grad g(\bx^k) - \bv^k\|^2  \\
	                & \quad + \frac{1}{2}\|\bv^k\|^2 \|\grad f(\bu^k)-\bs^k\|^2         - ( \frac{\lambda}{\tau_k} - C_fL_g - 2 )D_{h_g}(\bx^{k+1},\bx^k),
\end{align*}
where the second inequality follows from the Young's inequality and the last one uses the fact $h_g$ is 1-strongly convex.

\vspace{-0.5cm}
\subsection{Proof of Lemma \ref{lemma:RoR sequence bounded}}\label{pf:RoR sequence bounded}
By Assumption~\ref{assump:SoR stochastic oracle}, we have $\mE[\|\bu^k - g(\bx^k)\|^2 |\cF_k] \leq {\sigma_g^2}/{\abs{\cB_g^k}}$ and by Jensen's inequality, we get
$\mE[\|\bu^k - g(\bx^k)\| |\cF_k ] \leq \sqrt{\mE[\|\bu^k - g(\bx^k)\|^2 |\cF_k]} \leq {\sigma_g}/{\sqrt{|\cB_g^k|}}.$
Taking expectation on both sides of \eqref{eq:RoR ub} conditioned on $\cF_k$, we have
\begin{align*}
	 & \quad \mE[f(g(\bx^{k+1}))|\cF_k]  \\
	 & \leq f(g(\bx^k)) + 2C_f \mE[\|g(\bx^k) - \bu^k\||\cF_k]+ \frac{C_f^2}{2} \mE[\|\grad g(\bx^k) - \bv^k\|^2 |\cF_k]  \\
	 & \quad  +\frac{1}{2} \mE[\|\bv^k\|^2\|\grad f(\bu^k) - \bs^k\|^2|\cF_k ]  -\left( \frac{\lambda}{\tau} - C_fL_g -2 \right)\mE [D_{h_g}(\bx^{k+1},\bx^k) |\cF_k]  \\
	 & \leq f(g(\bx^k)) +  \frac{2C_f\sigma_g}{\sqrt{\abs{\cB_g^k}}} + \frac{C_f^2C_g^2}{2\abs{\cB_{\grad g}^k}} + \frac{C_f^2C_g^2}{2\abs{\cB_{\grad f}^k}}                                         - \left( \frac{\lambda}{\tau} - C_fL_g - 2 \right) \mE [D_{h_g}(\bx^{k+1},\bx^k) |\cF_k],
\end{align*}
where the second inequality holds because $\bv^k$, $\bu^k$, and $\bs^k$ are independent conditioned on $\cF_k$ and Remark~\ref{rem:stoch_bound_lip}.
Adding $-F^*$ on both sides and using the tower property, we get
\begin{align*}
	 & \quad  \mE[f(g(\bx^{k+1}))- F^*]  \\
	 & \leq \mE[f(g(\bx^k))-F^*] +  \frac{2C_f\sigma_g}{\sqrt{\abs{\cB_g^k}}} + \frac{C_f^2C_g^2}{2\abs{\cB_{\grad g}^k}} + \frac{C_f^2C_g^2}{2\abs{\cB_{\grad f}^k}}  -\left( \frac{\lambda}{\tau} - C_fL_g - 2 \right) \mE [D_{h_g}(\bx^{k+1},\bx^k) ].
\end{align*}
Rearranging the above terms, we have
\begin{align*}
	 & \quad  \left(\lambda \tau - (C_fL_g+2)\tau^2 \right)  \mE \left[{D_{h_g}(\bx^{k+1},\bx^k)}/{\tau^2}\right]  \\
	 & \leq -\mE[f(g(\bx^{k+1}))- F^*] +\mE[f(g(\bx^k))- F^*]  + \frac{2C_f\sigma_g}{\sqrt{\abs{\cB_g^k}}} + \frac{C_f^2C_g^2}{2\abs{\cB_{\grad g}^k}} + \frac{C_f^2C_g^2}{2\abs{\cB_{\grad f}^k}}  .\numberthis \label{eq:mini-batch error bound}
\end{align*}
Telescoping from $k=0$ to $K-1$,
{\small
\begin{align*}
	 \left(\lambda \tau - (C_fL_g+2)\tau^2 \right) \sum_{k=0}^{K-1} \mE \left[{D_{h_g}(\bx^{k+1},\bx^k)}/{\tau^2}\right]
	 \leq f(g(\bx^0)) - F^* + \sum_{k=0}^{K-1}\left( \frac{2C_f\sigma_g}{\sqrt{\abs{\cB_g^k}}} + \frac{C_f^2C_g^2}{2\abs{\cB_{\grad g}^k}} + \frac{C_f^2C_g^2}{2\abs{\cB_{\grad f}^k}} \right)
\end{align*}
}

Define $M_1 \triangleq \lambda \tau - (C_fL_g+2)\tau^2$ and let $\tau< \frac{C_fL_g+2C_{h_f}L_fL_g}{C_fL_g+2}$ which results in $M_1>0$. Setting $\abs{\cB_g^k} \equiv \left\lceil \frac{ 16C_f^2\sigma_g^2}{M_1^2\epsilon^2} \right\rceil$, $\abs{\cB_{\grad f}^k} = \abs{\cB_{\grad g}^k} \equiv \left\lceil \frac{2C_f^2C_g^2}{M_1\epsilon} \right\rceil $ and dividing both sides of the above inequality by $M_1K$, we get
\begin{align*}
	\frac{1}{K}\sum_{k=0}^{K-1} \mE \left[{D_{h_g}(\bx^{k+1},\bx^k)}/{\tau^2}\right] \leq \frac{f(g(\bx^0)) - F^*}{M_1K} + \epsilon.
\end{align*}

\vspace{-0.5cm}
\subsection{Proof of Lemma \ref{lemma:f-smooth} }\label{pf:f-smooth}
$h_f$ is twice differentiable and $L_{h_f}$-smooth over $\cS$, i.e.,
\begin{align*}
	\grad ^2h_f(\bx) \preceq L_{h_f}I,  \ \forall \bx \in \cS \subset \mathbb{R}^d,
\end{align*}
where $I$ is the identity matrix and $\cS$ is a bounded subset of $\mathbb{R}^d$. As $f$ is $L_{h_f}$-smooth relative to $h_f$, we have
\begin{align*}
	\grad f^2(\bx) \preceq L_f \grad ^2 h_f(\bx) \preceq L_fL_{h_f}I, \ \forall \bx \in \cS \subset \mathbb{R}^d,
\end{align*}
which means $f$ is $L_fL_{h_f}$-smooth  {over $\cS$}.

\vspace{-0.5cm}
\subsection{Proof of Lemma \ref{lemma:RoR-relationship of x_tilde and x}}\label{pf:RoR-relationship of x_tilde and x}
By the three-point inequality, $\forall \by \in \cX$, we have
\begin{align*}
	 & \tau_k\fprod{\bw^k,\bx^{k+1} - \by} + L_fh_f(\bu^k + (\bv^k)^\intercal(\bx^{k+1} - \bx^k)) - L_fh_f(\bu^k+(\bv^k)^\intercal(\by-\bx^k))  \\
	 & \quad + L_f \fprod{\bv^k\grad h_f(\bu^k),\by - \bx^{k+1}} \leq \lambda (D_{h_g}(\by,\bx^k) - D_{h_g}(\bx^{k+1},\bx^k) - D_{h_g}(\by,\bx^{k+1})    ),
\end{align*}
and
\begin{align*}
	 & \tau_k\fprod{\grad F(\bx^k),\tilde{\bx}^{k+1} - \by} + L_fh_f(g(\bx^k) + \grad g(\bx^k)^\intercal(\tilde{\bx}^{k+1} - \bx^k)) - L_fh_f(g(\bx^k)+\grad g(\bx^k)^\intercal(\by-\bx^k))  \\
	 & \quad + L_f \fprod{\grad g(\bx^k)\grad h_f(g(\bx^k)),\by - \tilde{\bx}^{k+1}} \leq \lambda (D_{h_g}(\by,\bx^k) - D_{h_g}(\tilde{\bx}^{k+1},\bx^k) - D_{h_g}(\by,\tilde{\bx}^{k+1})    ).
\end{align*}
Letting $\by = \tilde{\bx}^{k+1}$ in the first inequality, $\by = \bx^{k+1}$ in the second one, summing them together, and rearranging the terms, we have
\begin{align*} 
	     & \lambda(D_{h_g}(\tilde{\bx}^{k+1},\bx^{k+1}) + D_{h_g}(\bx^{k+1},\tilde{\bx}^{k+1}))  \\
	\leq & \tau_k \underbrace{ \langle\grad F(\bx^k) - \bw^k,\bx^{k+1} - \tilde{\bx}^{k+1}\rangle }_{T_1}+ L_f \underbrace{ \langle\bv^k\grad h_f(\bu^k) - \grad g(\bx^k)\grad h_f(g(\bx^k)),\bx^{k+1} - \tilde{\bx}^{k+1}\rangle}_{T_2} \\
	     & + L_f \underbrace{ \left[h_f(\bu^k + (\bv^k)^\intercal(\tilde{\bx}^{k+1}-\bx^k))  - h_f(g(\bx^k)+\grad g(\bx^k)^\intercal(\tilde{\bx}^{k+1} - \bx^k))\right]}_{T_3}  \\
	     & +L_f \underbrace{ \left[h_f(g(\bx^k)+\grad g(\bx^k)^\intercal(\bx^{k+1} - \bx^k)) -  h_f(\bu^k + (\bv^k)^\intercal(\bx^{k+1}-\bx^k))\right]}_{T_4}.\numberthis \label{eq:Le_Ri}
\end{align*}

We will next bound these four terms individually.
By Young's inequality,
\begin{align*}
	\tau_k T_1 & \leq \frac{\tau_k^2}{2C_fL_g}\|\bw^k - \grad F(\bx^k)\|^2 + \frac{C_fL_g}{2}\|\bx^{k+1} - \tilde{\bx}^{k+1}\|^2.
\end{align*}
Similarly,
\begin{align*}
	 & \quad L_f T_2  \\
	 & \leq \frac{L_f}{4C_{h_f}L_g}\|\bv^k\grad h_f(\bu^k) - \grad g(\bx^k) \grad h_f(g(\bx^k))\|^2 + C_{h_f}L_fL_g\|\bx^{k+1} - \tilde{\bx}^{k+1}\|^2  \\
	 & \quad + C_{h_f}L_fL_g\|\bx^{k+1} - \tilde{\bx}^{k+1}\|^2  \\
	        & \leq \frac{L_f}{2C_{h_f}L_g}\|\bv^k - \grad g(\bx^k)\|^2 \|\grad h_f(\bu^k) \|^2 + \frac{L_f}{2C_{h_f}L_g}\|\grad g(\bx^k)\|^2\|\grad h_f(\bu^k)  - \grad h_f(g(\bx^k)) \|^2 \\
	 & \quad + C_{h_f}L_fL_g\|\bx^{k+1} - \tilde{\bx}^{k+1}\|^2  \\
	        & \leq \frac{C_{h_f}L_f}{2L_g}\|\bv^k-\grad g(\bx^k)\|^2 + \frac{C_g^2L_fL_{h_f}^2}{2C_{h_f}L_g} \|\bu^k- g(\bx^k)\|^2 + C_{h_f}L_fL_g\|\bx^{k+1} - \tilde{\bx}^{k+1}\|^2.
\end{align*}
By Lipschitz continuity of $h_f$,
\begin{align*}
	L_f T_3 & \leq C_{h_f}L_f\|\bu^k - g(\bx^k) + (\bv^k-\grad g(\bx^k))^\intercal(\tilde{\bx}^{k+1}-\bx^k)\|  \\
	        & \leq C_{h_f}L_f \|\bu^k - g(\bx^k)\| + C_{h_f}L_f\|(\bv^k-\grad g(\bx^k))^\intercal(\tilde{\bx}^{k+1}-\bx^k)\|  \\
	        & \leq C_{h_f}L_f \|\bu^k - g(\bx^k)\| + \frac{C_{h_f}L_f}{2L_g}\|\bv^k-\grad g(\bx^k)\|^2 + \frac{C_{h_f}L_fL_g}{2}\|\tilde{\bx}^{k+1}-\bx^k\|^2.
	  \\
	L_f T_4 & \leq C_{h_f}L_f\| g(\bx^k)-\bu^k + (\grad g(\bx^k)-\bv^k)^\intercal(\bx^{k+1}-\bx^k)\|  \\
	        & \leq C_{h_f} L_f\|\bu^k - g(\bx^k)\| + C_{h_f}L_f\|(\bv^k-\grad g(\bx^k))^\intercal(\bx^{k+1}-\bx^k)\|  \\
	        & \leq C_{h_f}L_f \|\bu^k - g(\bx^k)\| + \frac{C_{h_f}L_f}{2L_g}\|\bv^k-\grad g(\bx^k)\|^2 + \frac{C_{h_f}L_fL_g}{2}\|\bx^{k+1}-\bx^k\|^2.
\end{align*}
$h_g$ is 1-strongly convex, so
\begin{align*}
	\lambda\|\tilde{\bx}^{k+1}-\bx^{k+1}\|^2 \leq \lambda (D_{h_g}(\tilde{\bx}^{k+1},\bx^{k+1}) + D_{h_g}(\bx^{k+1},\tilde{\bx}^{k+1})).
\end{align*}
Combining the upper and lower bounds we obtained for \eqref{eq:Le_Ri}, we have
\begin{align*}
	     & \quad     \lambda\|\tilde{\bx}^{k+1}-\bx^{k+1}\|^2  \\
	\leq & \frac{\tau_k^2}{2C_fL_g}\|\bw^k - \grad F(\bx^k)\|^2 + \frac{C_fL_g}{2}\|\bx^{k+1} - \tilde{\bx}^{k+1}\|^2  + \frac{C_{h_f}L_f}{2L_g}\|\bv^k - \grad g(\bx^k)\|^2 \\
	     & + \frac{C_g^2L_fL_{h_f}^2}{2C_{h_f}L_g}\|\bu^k- g(\bx^k)\|^2 + C_{h_f}L_fL_g\|\bx^{k+1} - \tilde{\bx}^{k+1}\|^2      + 2C_{h_f}L_f\|\bu^k-g(\bx^k)\|  \\
	     & + \frac{C_{h_f}L_f}{L_g}\norm{\bv^k-\grad g(\bx^k)}^2 + \frac{C_{h_f}L_fL_g}{2}\left(\|\tilde{\bx}^{k+1}-\bx^k\|^2 +\|\bx^{k+1}-\bx^k\|^2\right).
\end{align*}
Rearranging the terms and using $\lambda = C_fL_g+2C_{h_f}L_fL_g$, we have
\begin{align*}
	     & (\frac{C_fL_g}{2} +C_{h_f}L_fL_g ) \|\tilde{\bx}^{k+1}-\bx^{k+1}\|^2  \\
	\leq & \frac{\tau_k^2}{2C_fL_g}\|\bw^k - \grad F(\bx^k)\|^2 + \frac{3C_{h_f}L_f}{2L_g}\|\bv^k - \grad g(\bx^k)\|^2  + \frac{C_g^2L_fL_{h_f}^2}{2C_{h_f}L_g} \|\bu^k- g(\bx^k)\|^2 \\
	     & \quad + 2C_{h_f}L_f\|\bu^k-g(\bx^k)\| + \frac{C_{h_f}L_fL_g}{2}\left(\|\tilde{\bx}^{k+1}-\bx^k\|^2 +\|\bx^{k+1}-\bx^k\|^2  \right).
\end{align*}
Using $\norm{\bx-\by}^2 \leq 2 \norm{\bx-\bz}^2 + 2 \norm{\bz-\by}^2$ and above inequalities, we have
\begin{align*}
	 & \quad (\frac{C_fL_g}{2} +C_{h_f}L_fL_g )\|\tilde{\bx}^{k+1}-\bx^k\|^2  \\
	 & \leq 2(\frac{C_fL_g}{2} +C_{h_f}L_fL_g )\|\tilde{\bx}^{k+1}- \bx^{k+1}\|^2 + 2(\frac{C_fL_g}{2} +C_{h_f}L_fL_g )\|\bx^{k+1}-\bx^k\|^2  \\
	 & \leq \frac{\tau_k^2}{C_fL_g}\|\bw^k - \grad F(\bx^k)\|^2 + \frac{3C_{h_f}L_f}{L_g}\|\bv^k - \grad g(\bx^k)\|^2  + \frac{C_g^2L_fL_{h_f}^2}{C_{h_f}L_g} \|\bu^k- g(\bx^k)\|^2 \\
	 & \quad + 4C_{h_f}L_f\|\bu^k-g(\bx^k)\| + C_{h_f}L_fL_g\left(\|\tilde{\bx}^{k+1}-\bx^k\|^2 +\|\bx^{k+1}-\bx^k\|^2  \right)  \\
	 & \quad + 2(\frac{C_fL_g}{2} +C_{h_f}L_fL_g )\|\bx^{k+1}-\bx^k\|^2.
\end{align*}
Rearranging the terms,
	\begin{align*}
	 & \quad  \frac{C_fL_g}{2} \|\tilde{\bx}^{k+1}-\bx^k\|^2  \\
	 & \leq \frac{\tau_k^2}{C_fL_g}\|\bw^k - \grad F(\bx^k)\|^2  + 4C_{h_f}L_f\|\bu^k - g(\bx^k)\| + \frac{C_g^2L_fL_{h_f}^2}{C_{h_f}L_g}\|\bu^k - g(\bx^k)\|^2 \\
	 & \quad + \frac{3C_{h_f}L_f}{L_g}\|\bv^k - \grad g(\bx^k)\|^2 + \left(C_fL_g + 3C_{h_f}L_fL_g \right)\|\bx^{k+1}-\bx^k\|^2.
\end{align*}
Multiplying the inequality above by $2/(C_f L_g)$ and using 1-strong convexity of $h_g$, the proof is complete.

\vspace{-0.5cm}
\subsection{Proof of Theorem \ref{thm:RoR-sample complexity}}\label{pf:RoR-sample complexity}
Multiplying \eqref{eq:RoR ub} by $2(C_fL_g+3C_{h_f}L_fL_g)$, multiplying \eqref{eq:x_tilde-bound} by $C_f L_g M_1/2\tau^2$, and adding the two inequalities, we get
\begin{align*}
	 & \quad  \frac{C_fL_gM_1}{2} \frac{\|\tilde{\bx}^{k+1}-\bx^k\|^2}{\tau^2}  \\
	 & \leq (2C_fL_g+6C_{h_f}L_fL_g)(f(g(\bx^k))-f(g(\bx^{k+1})) ) + \frac{M_1}{C_fL_g}\|\bw^k - \grad F(\bx^k)\|^2  \\
	 & \quad +(\frac{4C_{h_f}L_fM_1}{\tau^2}+ 4C_f^2L_g + 12C_{h_f}C_fL_fL_g )\|\bu^k-g(\bx^k)\|                                                                                         + \frac{C_g^2L_fL_{h_f}^2M_1}{C_{h_f}L_g\tau^2}\|\bu^k - g(\bx^k)\|^2 \\
	 & \quad + (\frac{3C_{h_f}L_fM_1}{L_g\tau^2} + C_f^3L_g + 3C_{h_f}C_f^2L_fL_g)\|\bv^k - \grad g(\bx^k)\|^2  \\
	 & \quad + (C_fL_g + 3C{h_f}L_fL_g)\|\bv^k\|^2\|\grad f(\bu^k)-\bs^k\|^2. \numberthis \label{eq:x_tilde up}
\end{align*}
Taking expectation of both sides conditioned on $\cF^k$, we have
\begin{align*}
	 & \quad  \frac{C_fL_gM_1}{2} \frac{\|\tilde{\bx}^{k+1}-\bx^k\|^2}{\tau^2}  \\
	 & \leq (2C_fL_g+6C_{h_f}L_fL_g)(f(g(\bx^k))-\mE[f(g(\bx^{k+1}))|\cF^k] )  \\
	 & \quad +\frac{M_1}{C_fL_g}\left(\frac{2C_g^2L_f^2L_{h_f}^2\sigma_g^2}{\abs{\cB_g^k}} + \frac{2C_f^2C_g^2}{\abs{\cB^k_{\grad f}}} + \frac{2C_f^2C_g^2}{\abs{\cB^k_{\grad g}}}\right)  \\
	 & \quad + (\frac{4C_{h_f}L_fM_1}{\tau^2}+ 4C_f^2L_g + 12C_{h_f}C_fL_fL_g )\frac{\sigma_g}{\sqrt{\abs{\cB_g^k}}} + \frac{C_g^2L_fL_{h_f}^2M_1}{C_{h_f}L_g\tau^2}\frac{\sigma_g^2}{\abs{\cB_g^k}} \\
	 & \quad +(\frac{3C_{h_f}L_fM_1}{L_g\tau^2} + C_f^3L_g+ 3C_{h_f}C_f^2L_fL_g)\frac{C_g^2}{\abs{\cB_{\grad g}^k}}+(C_fL_g + 3C_{h_f}L_fL_g)\frac{C_f^2C_g^2}{\abs{\cB_{\grad f}^k}}.
\end{align*}
Taking expectation over the algorithm, applying the tower property, we get
\begin{align*}
	 & \quad\frac{C_fL_gM_1}{2}\mE\left[ \frac{\|\tilde{\bx}^{k+1}-\bx^k\|^2}{\tau^2}\right]  \\
	 & \leq (2C_fL_g+6C_{h_f}L_fL_g)\mE\left[f(g(\bx^k)) - f(g(\bx^{k+1}))\right] + \left(\frac{C_g^2L_f^2L_{h_f}^2M_1^3}{8C_f^3L_g} + \frac{C_g^2L_fL_{h_f}^2M_1^3}{16C_{h_f}C_f^2L_g\tau^2} \right)\epsilon^2 \\
	 & \quad + \left(\frac{2M_1^2}{C_fL_g} + \frac{C_{h_f}L_fM_1^2}{C_f\tau^2} + \frac{3C_{h_f}L_fM_1^2}{2C_f^2L_g\tau^2} + 2C_fL_gM_1 + 6C_{h_f}L_fL_gM_1 \right)\epsilon.
\end{align*}
Telescoping from $k=0$ to $K-1$, dividing both sides by $\frac{C_fL_gM_1K}{2}$, the result will follow.

\vspace{-0.5cm}
\subsection{Proof of Lemma \ref{lemma:RoS-sequence bound}}\label{pf:RoS-sequence bound}
By triangle inequality, we have
\begin{align*}
	\mE[\|\bv_j^k\|^2 |\cF_j^k] & \leq 2 \mE[\|\bv^k_j - \grad g(\bx_j^k)\|^2 |\cF_j^k] + 2 \mE[\|\grad g(\bx^k_j)\|^2|\cF_j^k] \\
	                                   & \leq 2 \mE[\|\bv^k_j - \grad g(\bx_j^k)\|^2 |\cF_j^k] + 2 C_g^2.\numberthis\label{eq:RoS-v bound}\end{align*}

Adjusting \eqref{eq:RoR ub} based on Algorithm \ref{alg:RoS-SPIDER}, we have
\begin{align*}
	f(g(\bx^k_{j+1})) & \leq f(g(\bx^k_j)) + 2C_f\|g(\bx^k_j)-\bu^k_j\| + \frac{C_f^2}{2}\|\grad g(\bx^k_j) - \bv^k_j\|^2 \\
	                  & \quad 
	                  + \frac{1}{2}\|\bv^k_j\|^2 \|\grad f(\bu^k_j)-\bs^k_j\|^2         -\frac{1}{2} ( \frac{\lambda}{\tau} - C_fL_g - 2 )\|\bx^k_{j+1}-\bx^k_j\|^2.
\end{align*}
Taking expectation of both sides, using the tower property and  \eqref{eq:RoS-v bound}, we have
\begin{align*}
	 & \quad \frac{1}{2} \mE[\|\bv^k_j\|^2\|\grad f(\bu^k_j) - \bs^k_j\|^2 ]  \\
	 & = \frac{1}{2}\mE\left[\mE[\|\bv^k_j\|^2|\cF^k_j]\mE[\|\grad f(\bu^k_j) - \bs^k_j\|^2 |\cF^k_j] \right]  \\
	 & \leq \frac{C_f^2}{2|\cB^{k,j}_{\grad f}|}\mE\left[(2\mE[\|\bv^k_j - \grad g(\bx^k_j)\|^2|\cF^k_j ] + 2C_g^2)  \right]  \\
	 & = \frac{C_f^2}{|\cB^{k,j}_{\grad f}|}\mE[ \|\bv^k_j - \grad g(\bx^k_j)\|^2] + \frac{C_f^2C_g^2}{|\cB^{k,j}_{\grad f}|}. \numberthis \label{eq:RoR-v times grad f}
\end{align*}
Hence, we have
{\small
\begin{align*}
	\mE[ f(g(\bx^k_{j+1}))] & \leq \mE[f(g(\bx^k_{j})) ] +  2C_f\mE[\|g(\bx^k_j) -\bu^k_j\|] + \left(\frac{C_f^2}{2} + \frac{C_f^2}{|{\cB^{k,j}_{\grad f}}|} \right)\mE[\|\grad g(\bx^k_j) - \bv^k_j\|^2] \\
	                               & \quad + \frac{C_f^2C_g^2}{\abs{\cB^{k,j}_{\grad f}}} -\frac{1}{2} \left( \frac{\lambda}{\tau} - C_fL_g - 2 \right)\mE[\|\bx^k_{j+1}-\bx^k_j\|^2].
\end{align*}
}

Applying Lemma \ref{lemma:RoS error} and inequality \eqref{eq:u-g_UB_spider}, we have
{\small
\begin{align*}
	 & \quad \mE[ f(g(\bx^k_{j+1}))]  \\
	 & \leq \mE[f(g(\bx^k_{j})) ] +  \frac{2C_f\sigma_g}{\sqrt{\abs{\cB^k_g}}} + C_f\delta  + \frac{C_f}{\delta}\sum_{r=0}^{j-1}\frac{C_g^2}{\abs{\cS^{k,r+1}_g}}\mE[\|\bx^k_{r+1} - \bx^k_r\|^2 ]   + \frac{C_f^2C_g^2}{2\abs{\cB^k_{\grad g}}} + \frac{C_f^2C_g^2}{\abs{\cB^{k,j}_{\grad f}}\abs{\cB^k_{\grad g}}} \\
	 & \quad + \left(\frac{C_f^2}{2} + \frac{C_f^2}{\abs{\cB^{k,j}_{\grad f}}} \right)\sum_{r=0}^{j-1}\frac{L_g^2}{\abs{\cS^{k,r+1}_{\grad g}}}\mE[\|\bx^k_{r+1} - \bx^k_r\|^2 ] + \frac{C_f^2C_g^2}{\abs{\cB^{k,j}_{\grad f}}}-\frac{1}{2} \left( \frac{\lambda}{\tau} - C_fL_g - 2 \right)\mE[\|\bx^k_{j+1}-\bx^k_j\|^2].
\end{align*}
}

Let $\abs{\cS_g^{k,j}}\equiv S_g, \abs{\cS_{\grad g}^{k,j}} \equiv S_{\grad g}, \abs{\cB^k_g}\equiv B_g, \abs{\cB_{\grad g}^k} \equiv B_{\grad g}, \abs{B^{k,j}_{\grad f}} \equiv B_{\grad f}$ and note that
\begin{align*}
	\sum_{r=0}^{j-1}\mE[\|\bx^k_{r+1}-\bx^k_r\|^2 ] \leq \sum_{r=0}^{J-1}\mE[\|\bx^k_{r+1}-\bx^k_r\|^2 ].
\end{align*}
Also, we have
\begin{align*}
	 & \quad \mE[ f(g(\bx^k_{j+1}))]  \\
	 & \leq \mE[f(g(\bx^k_{j})) ] + \left(\frac{C_fC_g^2}{\delta S_g} + \frac{C_f^2L_g^2}{2S_{\grad g}} + \frac{C_f^2L_g^2}{B_{\grad f}S_{\grad g}} \right)\sum_{j=0}^{J-1}\mE[\|\bx^k_{j+1}-\bx^k_j\|^2 ]  \\
	 & \quad + \frac{2C_f\sigma_g}{\sqrt{B_g}} + C_f\delta + \frac{C_f^2C_g^2}{2B_{\grad g}} + \frac{C_f^2C_g^2}{B_{\grad f}B_{\grad g}} + \frac{C_f^2C_g^2}{B_{\grad f}} -\frac{1}{2} \left( \frac{\lambda}{\tau} - C_fL_g - 2 \right)\mE[\|\bx^k_{j+1}-\bx^k_j\|^2].
\end{align*}
Telescoping the above inequality from $j=0$ to $J-1$ and rearranging the terms,
{\small
\begin{align*}
	 & \quad\left(\frac{\lambda}{2\tau} - \frac{C_fL_g}{2} - 1 - J\left(\frac{C_fC_g^2}{\delta S_g} + \frac{C_f^2L_g^2}{2S_{\grad g}} + \frac{C_f^2L_g^2}{B_{\grad f}S_{\grad g}} \right) \right)\sum_{j=0}^{J-1}\mE[\|\bx^k_{j+1}-\bx^k_j\|^2 ] \\
	 & \leq \mE[f(g(\bx^k_0)) - f(g(\bx^k_J))] + J\left(\frac{2C_f\sigma_g}{\sqrt{B_g}} + C_f\delta + \frac{C_f^2C_g^2}{2B_{\grad g}} + \frac{C_f^2C_g^2}{B_{\grad f}B_{\grad g}} + \frac{C_f^2C_g^2}{B_{\grad f}} \right).
\end{align*}
}
Telescoping  the above inequality from $k=0$ to $K-1$, and noting that by  Algorithm \ref{alg:RoS-SPIDER}, $\bx_0^{k+1} = \bx_J^k$, we have
{\small
\begin{align*}
	 & \quad\left(\frac{\lambda}{2\tau} - \frac{C_fL_g}{2} - 1 - J\left(\frac{C_fC_g^2}{\delta S_g} + \frac{C_f^2L_g^2}{2S_{\grad g}} + \frac{C_f^2L_g^2}{B_{\grad f}S_{\grad g}} \right) \right)\sum_{k=0}^{K-1}\sum_{j=0}^{J-1}\mE[\|\bx^k_{j+1}-\bx^k_j\|^2 ] \\
	 & \leq f(g(\bx^0_0)) - F^* + JK\left(\frac{2C_f\sigma_g}{\sqrt{B_g}} + C_f\delta + \frac{C_f^2C_g^2}{2B_{\grad g}} + \frac{C_f^2C_g^2}{B_{\grad f}B_{\grad g}} + \frac{C_f^2C_g^2}{B_{\grad f}} \right).
\end{align*}
}
Since $B_{\grad f} \geq 1$, the term $- J\left(\frac{C_fC_g^2}{\delta S_g} + \frac{C_f^2L_g^2}{2S_{\grad g}} + \frac{C_f^2L_g^2}{B_{\grad f}S_{\grad g}} \right) $ in the left hand side of the above inequality can be lower bounded as
{\small
\begin{align*}
	- J\left(\frac{C_fC_g^2}{\delta S_g} + \frac{C_f^2L_g^2}{2S_{\grad g}} + \frac{C_f^2L_g^2}{B_{\grad f}S_{\grad g}} \right) & \geq - J\left(\frac{C_fC_g^2}{\delta S_g} + \frac{3C_f^2L_g^2}{2S_{\grad g}} \right).
\end{align*}
}
Define $M_2 \triangleq \left(\frac{C_fL_g\tau}{2} - \left(\frac{C_fL_g}{2}+1\right)\tau^2\right)$, set  $\delta = \frac{M_2\epsilon}{3C_f}$, $S_g = \left \lceil\frac{6C_f^2C_g^2\tau}{C_{h_f}L_fL_g}\frac{J}{M_2\epsilon}\right\rceil$, $S_{\grad g} = \left\lceil \frac{3C_f^2L_g\tau J}{C_{h_f}L_f} \right \rceil$, $B_g = \left\lceil \frac{36C_f^2\sigma_g^2}{M_2^2\epsilon^2} \right \rceil$, $B_{\grad g} =B_{\grad f}=\left\lceil \frac{15C_f^2C_g^2}{2M_2\epsilon} \right \rceil$ and note that $\frac{C_f^2C_g^2}{B_{\grad f}B_{\grad g }}\leq \frac{C_f^2C_g^2}{B_{\grad f}}$. Hence, we have
{\small
\begin{align*}
	M_2\sum_{k=0}^{K-1}\sum_{j=0}^{J-1}\mE\left[\frac{\|\bx^k_{j+1}-\bx^k_j\|^2}{\tau^2} \right]\leq f(g(\bx^0_0)) -F^* + M_2JK \epsilon.
\end{align*}
}
With $\tau < \frac{C_fL_g}{C_fL_g+2}$, we have $M_2>0$. Dividing both sides by $M_2JK$, the proof is complete.

\vspace{-0.5cm}
\subsection{Proof of Theorem \ref{thm:RoS-sample complexity}}\label{pf:RoS-sample complexity}

By Lemma \ref{lemma:RoS-sequence bound}, for any given $K,J$ and $\epsilon$, the sequence $\{\bx^k\}$ generated by Algorithm \ref{alg:RoS-SPIDER} lies in a bounded set with probability 1. Hence Assumption \ref{assumption:hf-smooth} is applicable. Adjusting \eqref{eq:x_tilde up} based on Algorithm \ref{alg:RoS-SPIDER}, we have
{\small
\begin{align*}
	 & \quad  \frac{C_fL_gM_1}{2} \frac{\|\tilde{\bx}^k_{j+1}-\bx^k_j\|^2}{\tau^2}  \\
	 & \leq (2C_fL_g+6C_{h_f}L_fL_g )(f(g(\bx^k_j))-f(g(\bx^k_{j+1})) ) + \frac{M_1}{C_fL_g}\|\bw^k_j - \grad F(\bx^k_j)\|^2  \\
	 & \quad +(\frac{4C_{h_f}L_fM_1}{\tau^2}+ 4C_f^2L_g + 12C_{h_f}C_fL_fL_g )\|\bu^k_j-g(\bx^k_j)\|                                                                                           + \frac{C_g^2L_fL_{h_f}^2M_1}{C_{h_f}L_g\tau^2}\|\bu^k_j - g(\bx^k_j)\|^2 \\
	 & \quad  + (\frac{3C_{h_f}L_fM_1}{L_g\tau^2} + C_f^3L_g + 3C_{h_f}C_f^2L_fL_g )\|\bv^k_j - \grad g(\bx^k_j)\|^2  \\
	 & \quad + \left(C_fL_g + 3C_{h_f}L_fL_g \right)\|\bv^k_j\|^2\|\grad f(\bu^k_j)-\bs^k_j\|^2.
\end{align*}
}
Taking expectation of both sides, by Lemma \ref{lemma:RoS-w error}, inequality \eqref{eq:RoR-v times grad f}, and the fact $\abs{\cB_{\grad f}^{k,j}}\geq 1$,
\begin{align*}
	\frac{C_fL_gM_1}{2}\mE\left[ {\|\tilde{\bx}^k_{j+1}-\bx^k_j\|^2}/{\tau^2}   \right]
	 & \leq A_0 \mE[f(g(\bx^k_j))-f(g(\bx^k_{j+1})) ] + A_1\mE[\|\bu^k_j-g(\bx^k_j)\|]  \\
	 & \quad + A_2\mE[\|\bu^k_j-g(\bx^k_j)\|^2] + A_3\mE[ \|\bv^k_j - \grad g(\bx^k_j)\|^2] +   {A_4}/{|{\cB_{\grad f}^{k,j}}|},
\end{align*}
where \begin{align*}
	A_0 & \triangleq   2C_fL_g+6C_{h_f}L_fL_g ,\quad 	A_1 \triangleq \frac{4C_{h_f}L_fM_1}{\tau^2}+ 4C_f^2L_g + 12C_{h_f}C_fL_fL_g ,  \\
	A_2 & \triangleq  \frac{2C_g^2L_f^2L_{h_f}^2M_1}{C_fL_g}+ \frac{C_g^2L_fL_{h_f}^2M_1}{C_{h_f}L_g\tau^2}, \quad A_3 \triangleq   \frac{2C_fM_1}{L_g} + \frac{3C_{h_f}L_fM_1}{L_g\tau^2} + 3C_f^3L_g + 9C_{h_f}C_f^2L_fL_g, \\
	A_4 & \triangleq   \frac{2C_fC_g^2M_1}{L_g}+ 2C_f^3C_g^2L_g + 6C_{h_f}C_f^2C_g^2L_fL_g.
\end{align*}
By Lemma \ref{lemma:RoS error}, similar to the proof of Lemma \ref{lemma:RoS-sequence bound}, after telescoping twice, we get
{\small
\begin{align*}
	 & \quad\frac{1}{KJ}\sum_{k=0}^{K-1}\sum_{j=0}^{J-1} \mE\left[{\|\tilde{\bx}^k_{j+1}-\bx^k_j\|^2}/{\tau^2} \right]  \\
	 & \leq \frac{2A_0\left(f(g(\bx^0_0)) - F^* \right)}{C_fL_gM_1KJ} + \frac{2\tau^2}{C_fL_gM_1K}\left(\frac{A_1C_g^2}{2\delta S_g} + \frac{A_2C_g^2}{S_g} + \frac{A_3L_g^2}{S_{\grad g}} \right)\sum_{k=0}^{K-1}\sum_{j=0}^{J-1} \mE\left[ {\|\bx^k_{j+1} - \bx^k_j\|^2}/{\tau^2} \right] \\
	 & \quad + \frac{2}{C_fL_gM_1}\left(\frac{A_1\sigma_g}{\sqrt{B_g}}+ \frac{A_1\delta}{2} + \frac{A_2\sigma_g^2}{B_g} + \frac{A_3C_g^2}{B_{\grad g}} + \frac{A_4}{B_{\grad f}} \right).
\end{align*}
}
By Lemma \ref{lemma:RoS-sequence bound}, the right hand side of the above inequality can be further bounded as
\begin{align*}
	 & \quad  \frac{1}{KJ}\sum_{k=0}^{K-1}\sum_{j=0}^{J-1} \mE\left[{\|\tilde{\bx}^k_{j+1}-\bx^k_j\|^2}/{\tau^2} \right]  \\
	 & \leq \left(\frac{2A_0}{C_fL_g} + \frac{A_1C_g^2\tau^2J}{C_fL_g\delta S_gM_1}+ \frac{2A_2C_g^2\tau^2J}{C_fL_gS_gM_1} + \frac{2A_3L_g\tau^2J}{C_fS_{\grad g}M_1} \right)\frac{f(g(\bx^0_0)) - F^*}{M_1KJ} \\
	 & \quad + \left(  \frac{A_1C_g^2\tau^2J}{C_fL_g\delta S_gM_1}+ \frac{2A_2C_g^2\tau^2J}{C_fL_gS_gM_1} + \frac{2A_3L_g\tau^2J}{C_fS_{\grad g}M_1}\right)\epsilon  \\
	 & \quad + \frac{2}{C_fL_gM_1}\left(\frac{A_1\sigma_g}{\sqrt{B_g}}+ \frac{A_1\delta}{2} + \frac{A_2\sigma_g^2}{B_g} + \frac{A_3C_g^2}{B_{\grad g}} + \frac{A_4}{B_{\grad f}} \right).
\end{align*}
Replacing coefficients with the values set in Lemma \ref{lemma:RoS-sequence bound}, the final result is
\begin{align*}
	 & \quad\frac{1}{KJ}\sum_{k=0}^{K-1}\sum_{j=0}^{J-1} \mE\left[{\|\tilde{\bx}^k_{j+1}-\bx^k_j\|^2}/{\tau^2} \right]  \\
	 & \leq \left(\frac{2A_0}{C_fL_gM_1} + \frac{A_1C_{h_f}L_f\tau}{2C_f^2M_1^2} + \frac{2A_3C_{h_f}L_f\tau}{3C_f^3M_1^2} \right)\frac{f(g(\bx^0_0)) - F^*}{KJ}  \\
	 & \quad + \left(\frac{A_1C_{h_f}L_f\tau}{2C_f^2M_1} + \frac{2A_3C_{h_f}L_f\tau}{3C_f^3M_1} + \frac{2A_1M_2}{3C_f^2L_gM_1}  + \frac{4A_3M_2}{15C_f^3L_gM_1} + \frac{4A_4M_2}{15C_f^3C_g^2L_gM_1} \right)\epsilon \\
	 & \quad + \frac{A_2C_{h_f}L_fM_2\tau}{3C_f^3M_1^2} \frac{\left(f(g(\bx^0_0)) - F^*\right)\epsilon}{KJ}  + \left( \frac{A_2C_{h_f}L_fM_2\tau}{3C_f^3M_1} + \frac{A_2M_2^2}{18C_f^3L_gM_1} \right)\epsilon^2.
\end{align*}
To achieve $\epsilon$-stationarity, we need to set $KJ = \cO(\epsilon^{-1})$, then the number of calls to the $g_\xi$ oracle is
$$
	KB_g + KJS_g = \frac{\cO(\epsilon^{-1})}{J} \cO(\epsilon^{-2} ) + \cO(\epsilon^{-2})J,
$$
as $S_g = J\times \cO(\epsilon^{-1})$. Minimizing the right hand side of above equality w.r.t. $J$ results in  $J = \cO(\epsilon^{-1/2})$ and the sample complexity for $g_\xi$ equal to $\cO(\epsilon^{-5/2} )$.
Similarly, the total number of calls to the $\grad g_\xi$ oracle is
$$	
	KB_{\grad g} + KJS_{\grad g} = \frac{\cO(\epsilon^{-1})}{J}\cO(\epsilon^{-1}) + \cO(\epsilon^{-1})J,
$$
which achieves its minimum with $J = \cO(\epsilon^{-1/2})$ and the sample complexity of $\grad g_\xi$ is $\cO(\epsilon^{-3/2})$.

\vspace{-0.5cm}
\section{Experiment Proofs}
\subsection{Proof of Lemma \ref{lemma:4th rel-smooth}}
$f(\bx)$ is $L$-relative smooth to $h(\bx)$ is equivalent to $Lh(\bx) \pm f(\bx)$ is convex, i.e., $ L\grad^2 h(\bx) \pm \grad^2 f(\bx) \succeq 0$ since $h,f$ are twice differentiable. It is sufficient to show $\lambda_{\max}(\pm \grad^2 f(\bx)) \leq L \lambda_{\min}(\grad^2 h(\bx))$ where $\lambda_{\max}(A),\lambda_{\min}(A)$ are the maximal and minimum eigenvalues of matrix $A$. We have
$\grad^2 f(\bx) = 2A\bx \bx^\intercal A + \bx^\intercal A\bx A$ and $\grad^2h(\bx) = \|\bx\|_2^2I + 2 \bx \bx^\intercal$,
where $I$ is the identity matrix. Note that
$$\lambda_{\max}(\pm\grad^2f(\bx)) \leq \|\pm \grad^2f(\bx)\|_2 = \|\grad^2f(\bx)\|_2 \leq 3\|A\|_2^2\|\bx\|_2^2.$$
Let $L \geq 3\|A\|_2^2$, we have
$$ \lambda_{\max}(\pm\grad^2f(\bx)) \leq L\|\bx\|_2^2 \leq L \lambda_{\min}(\grad^2 h(\bx)),$$
which finishes the proof.

\vspace{-0.5cm}
\subsection{Proof of Lemma \ref{lem:rel_app_1}}

Let $g_1(\bx)\triangleq  \mE\left[\frac{1}{2}\bx^\intercal A_\xi \bx\right] $ and $g_2(\bx) \triangleq  \mE\left[\frac{1}{4}(\bx^\intercal A_\xi\bx)^2\right]$.
 {
Following the proof of Lemma \ref{lemma:4th rel-smooth}, we have 
\begin{align*}
	\grad^2 g_2(\bx) = \mE[2A_\xi\bx\bx^\intercal A_\xi + \bx^\intercal A_\xi \bx A_\xi].
\end{align*}
To  derive the sufficient condition $\lambda_{\max}(\pm \grad^2g_2(\bx))\leq L_g\lambda_{\min}(\grad^2 h(\bx))$, we have
\begin{align*}
		\lambda_{\max}(\pm \grad^2 g_2(\bx)) \leq \|\grad^2g_2(\bx)\|_2 = \|\mE[2A_\xi\bx\bx^\intercal A_\xi + \bx^\intercal A_\xi\bx A_\xi] \|_2 \leq 3\mE[\|A_\xi\|_2^2]\|\bx\|_2^2.
\end{align*}
	With $L \geq \mE[\|A_\xi\|_2^2]$, $h_2(\bx) = \frac{1}{4}\|\bx\|^4$, we can get
\begin{align*}
	\lambda_{\max}(\pm \grad^2 g_2(\bx)) \leq L\|\bx\|_2^2\leq L\lambda_{\min}(\grad ^2h_2(\bx)).
\end{align*}
}

 We have
\begin{align*}
		 & \quad   \|g(\bx) - g(\by)  - \grad g(\by)^\intercal(\bx-\by)    \|_2  \\
	 & \leq \|g(\bx) - g(\by)  - \grad g(\by)^\intercal(\bx-\by)    \|_1  \\
	 & =|g_1(\bx)-g_1(\by) - \fprod{\grad g_1(\by),\bx-\by}| + |g_2(\bx)-g_2(\by) - \fprod{\grad g_2(\by),\bx-\by}| \\
	 & \leq \frac{k_1}{2}\|\bx-\by\|_2^2 + k_2 D_{h_2}(\bx,\by)                                                 
	  =D_{h} (\bx,\by),
\end{align*}
	where $h(\bx)=\frac{k_1}{2}\norm{\bx}_2^2+\frac{k_2}{4}\norm{\bx}_2^4$, $k_1 \geq \norm{\mE\left[A_\xi\right]}$, $k_2 \geq 3\mE\left[\|A_\xi\|^2\right]$, and the second inequality is due to the smoothness of the quadric function  {and the relative smoothness of $g_2$}.


\end{appendices}

\end{document}